\documentclass[10.9pt,twoside]{amsart}
\usepackage{amsmath, amsthm, amscd, amsfonts, amssymb, graphicx, color}
\usepackage[bookmarksnumbered, colorlinks, plainpages]{hyperref}

\theoremstyle{definition}

\theoremstyle{remark}

\numberwithin{equation}{section}

\begin{document}
\setcounter{page}{1}
\begin{center}
{\bf  TOPOLOGICAL CENTERS  OF MODULE ACTIONS AND COHOMOLOGICAL GROUPS OF BANACH ALGEBRAS  }
\end{center}

\title[]{}
\author[]{KAZEM HAGHNEJAD AZAR   }

\address{}

\dedicatory{}



\begin{abstract} In this paper, first we  study some Arens regularity properties of module actions. Let  $B$ be a  Banach  $A-bimodule$ and let   ${Z}^\ell_{B^{**}}(A^{**})$ and  ${Z}^\ell_{A^{**}}(B^{**})$ be the topological centers of  the left module action $\pi_\ell:~A\times B\rightarrow B$ and the right module action $\pi_r:~B\times A\rightarrow B$, respectively. We investigate some relationships between topological center of $A^{**}$, ${Z}_1({A^{**}})$ with respect to the first Arens product and topological centers of module actions ${Z}^\ell_{B^{**}}(A^{**})$ and  ${Z}^\ell_{A^{**}}(B^{**})$. On the other hand, if $A$ has Mazure property and $B^{**}$ has the left $A^{**}-factorization$, then $Z^\ell_{A^{**}}(B^{**})=B$, and so for a locally compact non-compact group $G$ with compact covering number $card(G)$, we have $Z^\ell_{M(G)^{**}}{(L^1(G)^{**})}= {L^1(G)}$ and
$Z^\ell_{L^1(G)^{**}}{(M(G)^{**})}= {M(G)}$. By using the Arens regularity of module actions, we study some cohomological groups properties of Banach algebra and we extend some propositions from Dales,  Ghahramani,  Gr{\o}nb{\ae}k and others into general situations and we investigate the relationships between some cohomological groups of Banach algebra  $A$. We obtain some results in Connes-amenability of  Banach algebras, and so for  every compact group $G$, we conclude that  $H^1_{w^*}(L^\infty(G)^*,L^\infty(G)^{**})=0$.  Suppose that  $G$ is  an amenable locally compact group. Then there is a Banach $L^1(G)-bimodule$ such as $(L^\infty(G),.)$ such that $Z^1(L^1(G),L^\infty(G))=\{L_{f}:~f\in L^\infty(G)\}$ where for every $g\in L^1(G)$, we have $L_f(g)=f.g$.
\end{abstract} \maketitle

\section{\bf  Preliminaries and
Introduction }

\noindent As is well-known [1], the second dual $A^{**}$ of Banach algebra $A$ endowed with the either Arens multiplications is a Banach algebra. The constructions of the two Arens multiplications in $A^{**}$ lead us to definition of topological centers for $A^{**}$ with respect to both Arens multiplications. The topological centers of Banach algebras, module actions and applications of them  were introduced and discussed in [1, 5,  9, 14,  17, 18, 25, 26]. In first section, we  find some relationships between the topological centers of the second dual of Banach algebra $A$ and module actions with some conclusions in group algebras. For an unital Banach $A-module$ $B$ we show that ${Z}^\ell_{A^{**}}(B^{**}){Z}_1({A^{**}})={Z}^\ell_{A^{**}}(B^{**})$ and as results in group algebras, for locally compact group $G$,  we have ${Z}^\ell_{{L^1(G)}^{**}}(M(G)^{**})M(G)={Z}^\ell_{{L^1(G)}^{**}}(M(G)^{**})$ and ${Z}^\ell_{M(G)^{**}}({L^1(G)}^{**})M(G)={Z}^\ell_{M(G)^{**}}({L^1(G)}^{**})$. Let ${Z}^\ell_{A^{**}}(B^{**})A\subseteq B$ and suppose that $B$ is $WSC$, so we conclude that  ${Z}^\ell_{A^{**}}(B^{**})=B$. In [22], Neufang shows that if every Banach algebra $A$ satisfying $(M_k)$-property and whose dual $A^*$ has the property $(F_k)$, for some $k\geq \mathcal{N}_0$, then $A$ is left strongly Arens irregular and every linear left $A^{**}-module$ homomorphism on $A^*$ is automatically bounded and $weak^*-to-weak^*-continuous$. We extend this problem and some of its results into module actions and as its applications, we show that $Z^\ell_{M(G)^{**}}{(L^1(G)^{**})}= {L^1(G)}$ and
$Z^\ell_{L^1(G)^{**}}{(M(G)^{**})}= {M(G)}$ whenever $G$ is a locally compact non-compact group with covering number $card(G)$.\\
In second section, by using the Arens regularity of module actions, for Banach algebra $A$, we find some relations between cohomology groups $A$ and $A^{(2n)}$ with some applications in the $n-weak$ amenability of Banach algebras that introduced by Dales,  Ghahramani,   Gr{\o}nb{\ae}k in [6]. So for this aim, we extended some propositions from [6, 7, 10] into general situations. For  a Banach $A-bimodule$ $B$ and  $n\geq 0$, if the topological center of the left module action $\pi_\ell:A\times B\rightarrow B$ of $A^{(2n)}$ on $B^{(2n)}$ is $B^{(2n)}$  and  $H^1(A^{(2n+2)},B^{(2n+2)})=0$, then $H^1(A,B^{(2n)})=0$,  and we find the relationships between cohomological groups such as $H^1(A,B^{(n+2)})$  and    $H^1(A,B^{(n)})$,  spacial   $H^1(A,B^*)$ and  $H^1(A,B^{(2n+1)})$. We  investigated to relationships between cohomology groups of $A\oplus B$ and $A, B$ where $A$ and $B$ are Banach algebras In the section four, we establish some relationships between the Connes-amenability and supper-amenability of Banach algebras with some results in group algebras. In the section five,  by using the super-amenability of Banach algebra $A$, we give a representation  for the set of derivations from $A$ into $C$ where $C$ is a Banach $A-bimodule$. We have also some conclusions in the Arens regularity  of module actions. \\
We introduce some notations and definitions that we used
throughout  this paper.\\ Let $A$ be  a Banach algebra and $A^*$,
$A^{**}$, respectively, are the first and second dual of $A$.  For $a\in A$
 and $a^\prime\in A^*$, we denote by $a^\prime a$
 and $a a^\prime$ respectively, the functionals on $A^*$ defined by $\langle a^\prime a,b\rangle=\langle  a^\prime,ab\rangle=a^\prime(ab)$ and $\langle  a a^\prime,b\rangle=\langle  a^\prime,ba\rangle=a^\prime(ba)$ for all $b\in A$.
   The Banach algebra $A$ is embedded in its second dual via the identification
 $\langle  a,a^\prime\rangle$ - $\langle  a^\prime,a\rangle$ for every $a\in
A$ and $a^\prime\in
A^*$.  Let $A$ be a Banach algebra. We
say that a  net $(e_{\alpha})_{{\alpha}\in I}$ in $A$ is a left
approximate identity $(=LAI)$ [resp. right
approximate identity $(=RAI)$] if,
 for each $a\in A$,   $e_{\alpha}a\longrightarrow a$ [resp. $ae_{\alpha}\longrightarrow a$].\\
 \noindent Let $X,Y,Z$ be normed spaces and $m:X\times Y\rightarrow Z$ be a bounded bilinear mapping. Arens in [1] offers two natural extensions $m^{***}$ and $m^{t***t}$ of $m$ from $X^{**}\times Y^{**}$ into $Z^{**}$ as following\\
1. $m^*:Z^*\times X\rightarrow Y^*$,~~~~~given by~~~$\langle  m^*(z^\prime,x),y\rangle=\langle  z^\prime, m(x,y)\rangle$ ~where $x\in X$, $y\in Y$, $z^\prime\in Z^*$,\\
 2. $m^{**}:Y^{**}\times Z^{*}\rightarrow X^*$,~~given by $\langle  m^{**}(y^{\prime\prime},z^\prime),x\rangle=\langle  y^{\prime\prime},m^*(z^\prime,x)\rangle$ ~where $x\in X$, $y^{\prime\prime}\in Y^{**}$, $z^\prime\in Z^*$,\\
3. $m^{***}:X^{**}\times Y^{**}\rightarrow Z^{**}$,~ given by~ ~ ~$\langle  m^{***}(x^{\prime\prime},y^{\prime\prime}),z^\prime\rangle$  $=\langle  x^{\prime\prime},m^{**}(y^{\prime\prime},z^\prime)\rangle$\\ ~where ~$x^{\prime\prime}\in X^{**}$, $y^{\prime\prime}\in Y^{**}$, $z^\prime\in Z^*$.\\
The mapping $m^{***}$ is the unique extension of $m$ such that $x^{\prime\prime}\rightarrow m^{***}(x^{\prime\prime},y^{\prime\prime})$ from $X^{**}$ into $Z^{**}$ is $weak^*-to-weak^*$ continuous for every $y^{\prime\prime}\in Y^{**}$, but the mapping $y^{\prime\prime}\rightarrow m^{***}(x^{\prime\prime},y^{\prime\prime})$ is not in general $weak^*-to-weak^*$ continuous from $Y^{**}$ into $Z^{**}$ unless $x^{\prime\prime}\in X$. Hence the first topological center of $m$ may  be defined as following
$$Z_1(m)=\{x^{\prime\prime}\in X^{**}:~~y^{\prime\prime}\rightarrow m^{***}(x^{\prime\prime},y^{\prime\prime})~~is~~weak^*-to-weak^*-continuous\}.$$
Let now $m^t:Y\times X\rightarrow Z$ be the transpose of $m$ defined by $m^t(y,x)=m(x,y)$ for every $x\in X$ and $y\in Y$. Then $m^t$ is a continuous bilinear map from $Y\times X$ to $Z$, and so it may be extended as above to $m^{t***}:Y^{**}\times X^{**}\rightarrow Z^{**}$.
 The mapping $m^{t***t}:X^{**}\times Y^{**}\rightarrow Z^{**}$ in general is not equal to $m^{***}$, see [1], if $m^{***}=m^{t***t}$, then $m$ is called Arens regular. The mapping $y^{\prime\prime}\rightarrow m^{t***t}(x^{\prime\prime},y^{\prime\prime})$ is $weak^*-to-weak^*$ continuous for every $y^{\prime\prime}\in Y^{**}$, but the mapping $x^{\prime\prime}\rightarrow m^{t***t}(x^{\prime\prime},y^{\prime\prime})$ from $X^{**}$ into $Z^{**}$ is not in general  $weak^*-to-weak^*$ continuous for every $y^{\prime\prime}\in Y^{**}$. So we define the second topological center of $m$ as
$$Z_2(m)=\{y^{\prime\prime}\in Y^{**}:~~x^{\prime\prime}\rightarrow m^{t***t}(x^{\prime\prime},y^{\prime\prime})~~is~~weak^*-to-weak^*-continuous\}.$$
It is clear that $m$ is Arens regular if and only if $Z_1(m)=X^{**}$ or $Z_2(m)=Y^{**}$. Arens regularity of $m$ is equivalent to the following
$$\lim_i\lim_j\langle  z^\prime,m(x_i,y_j)\rangle=\lim_j\lim_i\langle  z^\prime,m(x_i,y_j)\rangle,$$
whenever both limits exist for all bounded sequences $(x_i)_i\subseteq X$ , $(y_i)_i\subseteq Y$ and $z^\prime\in Z^*$, see [5, 23].\\
 The regularity of a normed algebra $A$ is defined to be the regularity of its algebra multiplication when considered as a bilinear mapping. Let $a^{\prime\prime}$ and $b^{\prime\prime}$ be elements of $A^{**}$, the second dual of $A$. By $Goldstin^,s$ Theorem [4, P.424-425], there are nets $(a_{\alpha})_{\alpha}$ and $(b_{\beta})_{\beta}$ in $A$ such that $a^{\prime\prime}=weak^*-\lim_{\alpha}a_{\alpha}$ ~and~  $b^{\prime\prime}=weak^*-\lim_{\beta}b_{\beta}$. So it is easy to see that for all $a^\prime\in A^*$,
$$\lim_{\alpha}\lim_{\beta}\langle  a^\prime,m(a_{\alpha},b_{\beta})\rangle=\langle  a^{\prime\prime}b^{\prime\prime},a^\prime\rangle$$ and
$$\lim_{\beta}\lim_{\alpha}\langle  a^\prime,m(a_{\alpha},b_{\beta})\rangle=\langle  a^{\prime\prime}ob^{\prime\prime},a^\prime\rangle,$$
where $a^{\prime\prime}.b^{\prime\prime}$ and $a^{\prime\prime}ob^{\prime\prime}$ are the first and second Arens products of $A^{**}$, respectively, see [6, 18, 23].\\
The mapping $m$ is left strongly Arens irregular if $Z_1(m)=X$ and $m$ is right strongly Arens irregular if $Z_2(m)=Y$.\\
Regarding $A$ as a Banach $A-bimodule$, the operation $\pi:A\times A\rightarrow A$ extends to $\pi^{***}$ and $\pi^{t***t}$ defined on $A^{**}\times A^{**}$. These extensions are known, respectively, as the first (left) and the second (right) Arens products, and with each of them, the second dual space $A^{**}$ becomes a Banach algebra. In this situation, we shall also simplify our notations. So the first (left) Arens product of $a^{\prime\prime},b^{\prime\prime}\in A^{**}$ shall be simply indicated by $a^{\prime\prime}b^{\prime\prime}$ and defined by the three steps:
 $$\langle  a^\prime a,b\rangle=\langle  a^\prime ,ab\rangle,$$
  $$\langle  a^{\prime\prime} a^\prime,a\rangle=\langle  a^{\prime\prime}, a^\prime a\rangle,$$
  $$\langle  a^{\prime\prime}b^{\prime\prime},a^\prime\rangle=\langle  a^{\prime\prime},b^{\prime\prime}a^\prime\rangle.$$
 for every $a,b\in A$ and $a^\prime\in A^*$. Similarly, the second (right) Arens product of $a^{\prime\prime},b^{\prime\prime}\in A^{**}$ shall be  indicated by $a^{\prime\prime}ob^{\prime\prime}$ and defined by :
 $$\langle  a oa^\prime ,b\rangle=\langle  a^\prime ,ba\rangle,$$
  $$\langle  a^\prime oa^{\prime\prime} ,a\rangle=\langle  a^{\prime\prime},a oa^\prime \rangle,$$
  $$\langle a^{\prime\prime}ob^{\prime\prime},a^\prime\rangle=\langle b^{\prime\prime},a^\prime ob^{\prime\prime}\rangle.$$
  for all $a,b\in A$ and $a^\prime\in A^*$.\\
  The regularity of a normed algebra $A$ is defined to be the regularity of its algebra multiplication when considered as a bilinear mapping. Let $a^{\prime\prime}$ and $b^{\prime\prime}$ be elements of $A^{**}$, the second dual of $A$. By $Goldstine^,s$ Theorem [4, P.424-425], there are nets $(a_{\alpha})_{\alpha}$ and $(b_{\beta})_{\beta}$ in $A$ such that $a^{\prime\prime}=weak^*-\lim_{\alpha}a_{\alpha}$ ~and~  $b^{\prime\prime}=weak^*-\lim_{\beta}b_{\beta}$. So it is easy to see that for all $a^\prime\in A^*$,
$$\lim_{\alpha}\lim_{\beta}\langle a^\prime,\pi(a_{\alpha},b_{\beta})\rangle=\langle a^{\prime\prime}b^{\prime\prime},a^\prime\rangle$$ and
$$\lim_{\beta}\lim_{\alpha}\langle a^\prime,\pi(a_{\alpha},b_{\beta})\rangle=\langle a^{\prime\prime}ob^{\prime\prime},a^\prime\rangle,$$
where $a^{\prime\prime}b^{\prime\prime}$ and $a^{\prime\prime}ob^{\prime\prime}$ are the first and second Arens products of $A^{**}$, respectively, see [6, 18, 23].\\
We find the usual first and second topological center of $A^{**}$, which are
  $$Z_1(A^{**})=Z^\ell_1(A^{**})=\{a^{\prime\prime}\in A^{**}: b^{\prime\prime}\rightarrow a^{\prime\prime}b^{\prime\prime}~ is~weak^*-to-weak^*~continuous\},$$
   $$Z_2(A^{**})=Z_2^r(A^{**})=\{a^{\prime\prime}\in A^{**}: a^{\prime\prime}\rightarrow a^{\prime\prime}ob^{\prime\prime}~ is~weak^*-to-weak^*~continuous\}.$$\\
 An element $e^{\prime\prime}$ of $A^{**}$ is said to be a mixed unit if $e^{\prime\prime}$ is a
right unit for the first Arens multiplication and a left unit for
the second Arens multiplication. That is, $e^{\prime\prime}$ is a mixed unit if
and only if,
for each $a^{\prime\prime}\in A^{**}$, $a^{\prime\prime}e^{\prime\prime}=e^{\prime\prime}o a^{\prime\prime}=a^{\prime\prime}$. By
[4, p.146], an element $e^{\prime\prime}$ of $A^{**}$  is  mixed
      unit if and only if it is a $weak^*$ cluster point of some BAI $(e_\alpha)_{\alpha \in I}$  in
      $A$.\\
Let now $B$ be a Banach $A-bimodule$, and let
$$\pi_\ell:~A\times B\rightarrow B~~~and~~~\pi_r:~B\times A\rightarrow B.$$
be the right and left module actions of $A$ on $B$. Then $B^{**}$ is a Banach $A^{**}-bimodule$ with the following module actions where $A^{**}$ is equipped with the left Arens product
$$\pi_\ell^{***}:~A^{**}\times B^{**}\rightarrow B^{**}~~~and~~~\pi_r^{***}:~B^{**}\times A^{**}\rightarrow B^{**}.$$
Similarly, $B^{**}$ is a Banach $A^{**}-bimodule$ with the following module actions where $A^{**}$ is equipped with the right Arens product\\
$$\pi_\ell^{t***t}:~A^{**}\times B^{**}\rightarrow B^{**}~~~and~~~\pi_r^{t***t}:~B^{**}\times A^{**}\rightarrow B^{**}.$$\\
Let $B$ be a   Banach $A-bimodule$.
   A derivation from $A$ into $B$ is a bounded linear mapping $D:A\rightarrow B$ such that $$D(xy)=xD(y)+D(x)y~~for~~all~~x,~y\in A.$$
The space of continuous derivations from $A$ into $B$ is denoted by $Z^1(A,B)$.\\
Easy example of derivations are the inner derivations, which are given for each $b\in B$ by
$$\delta_b(a)=ab-ba~~for~~all~~a\in A.$$
The space of inner derivations from $A$ into $B$ is denoted by $N^1(A,B)$.
The Banach algebra $A$ is said to be a amenable, when for every Banach $A-bimodule$ $B$, the inner derivations are only derivations existing from $A$ into $B^*$. It is clear that $A$ is amenable if and only if $H^1(A,B^*)=Z^1(A,B^*)/ N^1(A,B^*)=\{0\}$. The concept of amenability for a Banach algebra $A$, introduced by Johnson in 1972, has proved to be of enormous importance in Banach algebra theory, see [14].
A Banach algebra $A$ is said to be a weakly amenable, if every derivation from $A$ into $A^*$ is inner. Similarly, $A$ is weakly amenable if and only if $H^1(A,A^*)=Z^1(A,A^*)/ N^1(A,A^*)=\{0\}$. The concept of weak amenability was first introduced by Bade, Curtis and Dales in [2] for commutative Banach algebras, and was extended to the noncommutative case by Johnson in [15].\\
For Banach $A-bimodule$ $B$, the quotient space $H^1(A,B)$ of all continuous derivations from $A$ into $B$ modulo the subspace of inner derivations is called the first cohomology group of $A$ with coefficients in $B$.\\

\begin{center}
\section{ \bf Topological centers of module actions }
\end{center}
In this section, the notations $WSC$ is used for weakly sequentially
complete Banach space $A$, that is, $A$ is said to be weakly
sequentially complete, if every weakly Cauchy sequence in $A$ has
a weak limit in $A$.\\
\noindent Suppose that $A$ is a Banach algebra and $B$ is a Banach $A-bimodule$. According to [6] $B^{**}$ is a Banach $A^{**}-bimodule$, where  $A^{**}$ is equipped with the first Arens product.
 We define  $B^{*}B$ as a subspace of $A$, that is, for all $b^{\prime}\in B^{*}$ and $b\in B$, we define
$$\langle b^{\prime}b,a\rangle =\langle b^{{\prime}},ba\rangle ;$$
We similarly define $B^{***}B^{**}$ as a subspace of $A^{**}$  and we take $A^{(0)}=A$ and $B^{(0)}=B$.\\
In the following, we will study some properties of topological centers of module actions and we will extend some problems from topological centers of Banach algebras into module actions with some relationships of them. We have some results in group algebras, that is, for  compact $G$ we have

$$Z_{L^1(G)^{**}}(M(G)^{**})L^1(G)\subseteq M(G)^{**}L^1(G)\subseteq L^1(G)=L^1(G)\subseteq M(G),$$
$$Z_{L^1(G)^{**}}(M(G)^{**})=M(G),$$
and for locally compact $G$ we obtain the following results:

$${Z}^\ell_{{L^1(G)}^{**}}(M(G)^{**})M(G)={Z}^\ell_{{L^1(G)}^{**}}(M(G)^{**}),$$

$${Z}^\ell_{M(G)^{**}}({L^1(G)}^{**})M(G)={Z}^\ell_{M(G)^{**}}({L^1(G)}^{**}),$$
$$Z_{M(G)^{**}}{(L^1(G)^{**})}=L^1(G).$$
\\

\noindent If we set $\pi_\ell (a,b)=ab$ and $\pi_r(b,a)=ba$, then we use the notions ${Z}^\ell_{B^{**}}(A^{**})$ and ${Z}^\ell_{A^{**}}(B^{**})$ for topological centers of
 module actions as follows.

$${Z}^\ell_{B^{**}}(A^{**})=\{a^{\prime\prime}\in A^{**}:~the~map~~b^{\prime\prime}\rightarrow a^{\prime\prime} b^{\prime\prime}~:~B^{**}\rightarrow B^{**}$$$$~is~~~weak^*-to-weak^*~continuous\}.$$
$${Z}^\ell_{A^{**}}(B^{**})=\{b^{\prime\prime}\in B^{**}:~the~map~~a^{\prime\prime}\rightarrow b^{\prime\prime} a^{\prime\prime}~:~A^{**}\rightarrow B^{**}$$$$~is~~~weak^*-to-weak^*~continuous\}$$\\

Lau and  \"{U}lger in [18], have studied some properties of the topological centers of second dual of Banach algebra $A$ and $(A^*A)^*$. They established some relations between them and subspaces of $A^{**}$ and $A^*$ with some results in group algebras. In the following, we extend some propositions from [18] into module actions with some new conclusions in group algebras.\\\\

\noindent {\it{\bf Theorem 2-1.}} We have the following assertions:
\begin{enumerate}
\item Let  $B$ be a  Banach left $A-module$ and $X=\{a^{\prime\prime}\in A^{**}:~B^*a^{\prime\prime}\subseteq B^*\}$. Then  $~{Z}^\ell_{B^{**}}(A^{**})X\subseteq {Z}^\ell_{B^{**}}(A^{**})$.
\item Let  $B$ be a  Banach right $A-module$ and  ${Z}^\ell_{A^{**}}(B^{**}){Z}_1({A^{**}})\subseteq {Z}^\ell_{A^{**}}(B^{**})$. If $B$ is an unital Banach $A-module$, then ${Z}^\ell_{A^{**}}(B^{**}){Z}_1({A^{**}})={Z}^\ell_{A^{**}}(B^{**})$.
\item  Let $A$ be a Banach algebra and $B\subseteq A$ be a Banach space such that $A$ is a Banach left $A-module$. If $A$ is Arens regular, then $~{Z}^\ell_{B^{**}}(A^{**})=A^{**}$.
\item    Suppose that $B$ is a Banach space including of a Banach algebra $A$ such that $B$ is a Banach left $A-module$. If $A$ is strongly Arens irregular, then $Z^\ell_{B^{**}}(A^{**})=A$. \\
    
\end{enumerate}

\begin{proof}
\begin{enumerate}
 \item ~Let  $a^{\prime\prime}\in {Z}^\ell_{B^{**}}(A^{**})$ and $x^{\prime\prime}\in X$. Suppose that $(b_{\alpha}^{\prime\prime})_{\alpha}\subseteq B^{**}$ such that  $b^{\prime\prime}_{\alpha} \stackrel{w^*} {\rightarrow}b^{\prime\prime}$. Since $b^\prime x^{\prime\prime}\in B^*$ for each $b^\prime\in B^*$, we have
$$\langle x^{\prime\prime}b_{\alpha}^{\prime\prime}, b^\prime\rangle =\langle b_{\alpha}^{\prime\prime}, b^\prime x^{\prime\prime}\rangle
\rightarrow \langle b^{\prime\prime}, b^\prime x^{\prime\prime}\rangle =\langle x^{\prime\prime}b^{\prime\prime}, b^\prime\rangle .$$
It follows that $x^{\prime\prime}b^{\prime\prime}_{\alpha} \stackrel{w^*} {\rightarrow}x^{\prime\prime}b^{\prime\prime}$. Then we conclude that
$$\langle (a^{\prime\prime}x^{\prime\prime})b_{\alpha}^{\prime\prime}, b^\prime\rangle =\langle a^{\prime\prime}(x^{\prime\prime}b_{\alpha}^{\prime\prime}), b^\prime\rangle \rightarrow
\langle a^{\prime\prime}(x^{\prime\prime}b^{\prime\prime}), b^\prime\rangle =\langle (a^{\prime\prime}x^{\prime\prime})b^{\prime\prime}, b^\prime\rangle .$$
It follows that $a^{\prime\prime}x^{\prime\prime}\in {Z}^\ell_{B^{**}}(A^{**})$.\\

\item ~~Let  $b^{\prime\prime}\in {Z}^\ell_{A^{**}}(B^{**})$  and $a^{\prime\prime}\in {Z}_1(A^{**})$. Assume that
$(x_{\alpha}^{\prime\prime})_{\alpha}\subseteq A^{**}$ such that  $x^{\prime\prime}_{\alpha} \stackrel{w^*} {\rightarrow}x^{\prime\prime}$. By using [18, Lemma 3.1], we have $a^\prime a^{\prime\prime}\in A^*$ for all $a^\prime \in A^*$. Then
$$\langle a^{\prime\prime}x^{\prime\prime}_\alpha,a^\prime\rangle =\langle x^{\prime\prime}_\alpha,a^\prime a^{\prime\prime}\rangle \rightarrow \langle x^{\prime\prime},a^\prime a^{\prime\prime}\rangle =\langle a^{\prime\prime}x^{\prime\prime},a^\prime \rangle .$$
It follows that $a^{\prime\prime}x^{\prime\prime}_{\alpha} \stackrel{w^*} {\rightarrow}a^{\prime\prime}x^{\prime\prime}$. Consequently we have $(b^{\prime\prime}a^{\prime\prime})x^{\prime\prime}_{\alpha} \stackrel{w^*} {\rightarrow}(b^{\prime\prime}a^{\prime\prime})x^{\prime\prime}$, and so $b^{\prime\prime}a^{\prime\prime}\in {Z}^\ell_{A^{**}}(B^{**})$. It follows that  ${Z}^\ell_{A^{**}}(B^{**}){Z}_1({A^{**}})\subseteq {Z}^\ell_{A^{**}}(B^{**})$.\\
Now let  $b^{\prime\prime}\in {Z}^\ell_{A^{**}}(B^{**})$ and suppose that $e\in A$ such that $eb=be=b$ for all $b\in B$.  Then for every $b^\prime\in B^*$ and $b\in B$, we have
$$\langle b^\prime e,b\rangle =\langle b^\prime ,eb\rangle =\langle b^\prime ,b\rangle .$$
It follows that $b^\prime e=b^\prime$. It is similar $eb^\prime =b^\prime$, and so
$$\langle b^{\prime\prime} e,b^\prime\rangle =\langle b^{\prime\prime}, eb^\prime\rangle =\langle b^{\prime\prime}, b^\prime\rangle .$$
Since $A\subseteq Z_1(A^{**})$, $e\in Z_1(A^{**})$.
Consequently $b^{\prime\prime} =b^{\prime\prime} e\in {Z}^\ell_{A^{**}}(B^{**}){Z}_1({A^{**}})$.
Hence ${Z}^\ell_{A^{**}}(B^{**}){Z}_1({A^{**}})={Z}^\ell_{A^{**}}(B^{**})$.
\item Since $A$ is Arens regular, $Z_1(A^{**})=A^{**}$. Therefore we have $A^{**}=Z_1(A^{**})\subseteq{Z}^\ell_{B^{**}}(A^{**})\subseteq A^{**}.$

\item  We know that $Z^\ell_{B^{**}}(A^{**})\subseteq Z_{1}(A^{**})$ and $A\subseteq Z^\ell_{B^{**}}(A^{**})$, so $Z^\ell_{B^{**}}(A^{**})=A$.
 
\end{enumerate}
\end{proof}

\noindent {\it{\bf Corollary 2-2.}} Let  $B$ be a  Banach  $A-bimodule$. Then we have the following assertions
\begin{enumerate}
\item $~{Z}^\ell_{B^{**}}(A^{**})A\subseteq {Z}^\ell_{B^{**}}(A^{**})$, and if $A$ is unital, then $$~{Z}^\ell_{B^{**}}(A^{**})A={Z}^\ell_{B^{**}}(A^{**}){Z}_1({A^{**}})={Z}^\ell_{B^{**}}(A^{**}).$$
\item  ${Z}^\ell_{A^{**}}(B^{**})A\subseteq {Z}^\ell_{A^{**}}(B^{**})$, and if $B$ is an unital as Banach $A-module$, then $${Z}^\ell_{A^{**}}(B^{**})A={Z}^\ell_{A^{**}}(B^{**}){Z}_1({A^{**}})={Z}^\ell_{A^{**}}(B^{**}).$$\\
\end{enumerate}

\noindent {\it{\bf Example 2-3.}} Let $G$ be a locally compact group. Then 
\begin{enumerate}
\item By using [17], we know that $Z_1({L^1(G)}^{**})= L^1(G)$ and $Z_1(M(G)^{**})= M(G)$. Since $L^1(G)$ is a Banach $M(G)$-bimodule and $M(G)$ is unital, by using Corollary 2-2, we have
$${Z}^\ell_{{L^1(G)}^{**}}(M(G)^{**})M(G)={Z}^\ell_{{L^1(G)}^{**}}(M(G)^{**}).$$
\item Since $M(G)$ also is  a Banach $L^1(G)$-bimodule and $M(G)$ is a unital, by using Corollary 2-2, we have
$${Z}^\ell_{M(G)^{**}}({L^1(G)}^{**})M(G)={Z}^\ell_{M(G)^{**}}({L^1(G)}^{**}).$$
\item  For locally compact group $G$, since $Z^\ell_{L^1(G)^{**}}{(L^1(G)^{**})}= {L^1(G)^{}}$ and  $L^1(G)$ is a subspace of $M(G)$, by using Theorem 2-1, it is clear that
$$Z^\ell_{M(G)^{**}}{(L^1(G)^{**})}=L^1(G).$$
\item    Consider the algebra $c_0=(c_0,.)$ is the collection of all sequences of scalars that convergence to $0$, with the some vector space operations and norm as $\ell^\infty$. By using [5, Example 2.6.22 (iii)] we know that $\ell^1=(\ell^1,.)$ is Arens regular, so by  Theorem 2-1, it is clear that $$Z^\ell_{c_0(G)^{**}}{(\ell^1(G)^{**})}=\ell^1(G)^{**}.$$\\

\end{enumerate}

\noindent {\it{\bf Theorem 2-4.}} Let  $B$ be a  Banach  $A-bimodule$. Then we have the following assertions:
\begin{enumerate}
\item $B{Z}_1(A^{**})\subseteq {Z}^\ell_{A^{**}}(B^{**})$.
\item  Let  $B^*B^{**}\subseteq A^*$. Then $~B^{**}{Z}_1(A^{**})\subseteq {Z}^\ell_{A^{**}}(B^{**})$. Moreover if  $B$ is an unital Banach $A-module$, then $B^{**}{Z}_1({A^{**}})={Z}^\ell_{A^{**}}(B^{**})$.
\end{enumerate}

\begin{proof}
\begin{enumerate}
\item  Let $b\in B$ and $a^{\prime\prime}\in Z_1(A^{**})$. Assume that
$(x_{\alpha}^{\prime\prime})_{\alpha}\subseteq A^{**}$ such that  $x^{\prime\prime}_{\alpha} \stackrel{w^*} {\rightarrow}x^{\prime\prime}$. Then for every $b^\prime\in B^*$, we have
$$\langle (ba^{\prime\prime})x_{\alpha}^{\prime\prime},b^\prime\rangle =\langle b(a^{\prime\prime}x_{\alpha}^{\prime\prime}),b^\prime\rangle =
\langle a^{\prime\prime}x_{\alpha}^{\prime\prime},b^\prime b\rangle \rightarrow \langle a^{\prime\prime}x^{\prime\prime},b^\prime b\rangle $$$$=\langle (ba^{\prime\prime})x^{\prime\prime},b^\prime \rangle .$$
It follows that $ba^{\prime\prime}\in {Z}^\ell_{A^{**}}(B^{**})$.
\item Suppose that $b^{\prime\prime}\in B^{**}$ and  $a^{\prime\prime}\in Z_1(A^{**})$. Let
$(x_{\alpha}^{\prime\prime})_{\alpha}\subseteq A^{**}$ such that  $x^{\prime\prime}_{\alpha} \stackrel{w^*} {\rightarrow}x^{\prime\prime}$. Since for every $b^\prime\in B^*$, $b^\prime b^{\prime\prime}\in A^{*}$, we have
$$\langle (b^{\prime\prime}a^{\prime\prime})x_{\alpha}^{\prime\prime},b^\prime\rangle
=\langle b^{\prime\prime}(a^{\prime\prime}x_{\alpha}^{\prime\prime}),b^\prime\rangle =\langle a^{\prime\prime}x_{\alpha}^{\prime\prime},b^\prime b^{\prime\prime}\rangle\rightarrow
\langle a^{\prime\prime}x^{\prime\prime},b^\prime b^{\prime\prime}\rangle =\langle (b^{\prime\prime}a^{\prime\prime})x^{\prime\prime},b^\prime\rangle .$$
It follows that $b^{\prime\prime}a^{\prime\prime}\in {Z}^\ell_{A^{**}}(B^{**})$.\\
Next part is clear.\end{enumerate}\end{proof}

\noindent {\it{\bf Example 2-5.}} Let $G$ be a locally compact group. Then \\
i)   For $1\leq p\leq\infty$, set $A=M(G)$ and $B=L^p(G)$. Thus by using proceeding theorem, we have
$$L^p(G)*M(G)=L^p(G)*Z_1(M(G)^{**})\subseteq Z^\ell_{M(G)^{**}}(L^p(G)).$$
ii) Take $A=L^1(G)$ and $B=L^\infty (G)$. Therefore by using proceeding theorem, we have
$$LUC(G)=L^\infty (G)*L^1(G)=L^\infty (G)*Z_1(L^1(G)^{**})\subseteq Z^\ell_{L^1(G)^{**}}(L^\infty(G)^{**}).$$
iii) In the proceeding theorem, if we take $B=c_0(G)$ and $A=\ell^1(G)$, then it is clear that $B$ is a Banach $A-bimodule$. So, by using [5, Example 2.6.22 (iii)]], since $(\ell^1(G),.)$ is Arens regular, we have the following relations

$c_0(G).\ell^1(G)\subseteq Z^\ell_{(\ell^\infty(G))^*}(\ell^\infty(G)).$
On the other hand, since $$\ell_1(G)\ell^\infty(G)=c_0(G)^*c_0(G)^{**}\subseteq \ell^\infty(G),$$ we have
$$\ell^\infty(G)\ell^1(G)=\ell^\infty(G) Z_1((\ell^1)^{**})\subseteq Z^\ell_{(\ell^1(G))^{**}}(\ell^\infty(G)).$$\\

\noindent {\it{\bf Theorem 2-6.}} Suppose that $B$ is a Banach $A-bimodule$ and $A$ has a sequential $BAI$ $(e_n)_{n\in \mathbb{N}}\subseteq A$. Then we have the following assertions.
\begin{enumerate}
\item  Let $AB^{**}\subseteq B$ and $B^*$ left factors with respect to $A$. If $B$ is $WSC$, then $B$ is reflexive.
\item Let ${Z}^\ell_{A^{**}}(B^{**})A\subseteq B$. If $B$ is $WSC$, then ${Z}^\ell_{A^{**}}(B^{**})=B$.\\
\end{enumerate}
\begin{proof}
\begin{enumerate}
\item  Let $b^{\prime\prime}\in B^{**}$ and suppose that $e^{\prime\prime}\in A^{**}$ is a mixed unit for $A^{**}$ such that it is $weak^*$ cluster of the sequence $(e_n)_{n\in \mathbb{N}}\subseteq A$. Without loss generality, we take $w^*-\lim_n e_n=e^{\prime\prime}$. Since $B^*$ left factors with respect to $A$, for every $b^\prime\in B^*$, there are $x^\prime\in B^*$ and $a\in A$ such that $b^\prime=x^\prime a$. Then we have
    $$\lim_n\langle e_nb^{\prime\prime},b^\prime\rangle =\lim_n\langle b^{\prime\prime},b^\prime e_n\rangle =\lim_n\langle b^{\prime\prime},x^\prime a e_n\rangle =\lim_n\langle b^{\prime\prime}x^\prime, a e_n\rangle $$$$\rightarrow \langle b^{\prime\prime}x^\prime, a \rangle =\langle b^{\prime\prime},x^\prime a \rangle =\langle b^{\prime\prime},b^\prime \rangle .$$
It follows that $w^*-\lim_ne_nb^{\prime\prime}=b^{\prime\prime}$ in $B^{**}$. Since $e_nb^{\prime\prime}\in B$ and $B$ is $WSC$, $b^{\prime\prime}\in B$, and so $B$ is reflexive.

\item Let $b^{\prime\prime}\in {Z}^\ell_{A^{**}}(B^{**})$. Then $w^*-\lim_n b^{\prime\prime}e_n=b^{\prime\prime}$  in $B^{**}$. Since ${Z}^\ell_{A^{**}}(B^{**})A\subseteq B$ and $B$ is $WSC$, $b^{\prime\prime}\in B$. Thus ${Z}^\ell_{A^{**}}(B^{**})=B$.\\
\end{enumerate}
\end{proof}

 In  part (2) of Theorem 2-6, if we take $B=A$, we obtain Theorem 3.4 (a) from [18].\\\\

 Lau and Losert in [17], for locally compact group $G$, show that the topological center of $L^1(G)^{**}$ is $L^1(G)$. Neufang in  [21] and [22] has studied the topological centers of $L^1(G)^{**}$ and $M(G)^{**}$ in spacial case, Eshaghi and Filali in [9] have studied these problems on module actions. In the following example for a compact group $G$, by using the proceeding theorem, we study the topological centers of $L^1(G)^{**}$, $\pi_\ell^{***}: L^1(G)^{**}\times M(G)^{**}\rightarrow L^1(G)^{**}$ and $\pi_r^{***}: M(G)^{**}\times L^1(G)^{**}\rightarrow M(G)^{**}$ with some new results.\\\\

\noindent {\it{\bf Example 2-7.}} Let $G$ be a  compact group. It is clear that a group algebra $L^1(G)$ has a sequential $BAI$. We know that $L^1(G)$ is an ideal in its second dual  $L^1(G)^{**}$ and it is $WSC$ Banach algebra. Thus we have
$$Z_1(L^1(G)^{**})L^1(G)\subseteq L^1(G)^{**}L^1(G)\subseteq L^1(G).$$
Then by using Theorem 2-6 part (2), we conclude that $Z_1(L^1(G)^{**})=L^1(G)$.\\
Now, if we take $A=L^1(G)$ and $B=M(G)$, then we have
$$Z^\ell_{L^1(G)^{**}}(M(G)^{**})L^1(G)\subseteq M(G)^{**}L^1(G)\subseteq L^1(G)=L^1(G)\subseteq M(G).$$
Since $M(G)$ is a $WSC$, by using Theorem 2-6 (2), we have $Z^\ell_{L^1(G)^{**}}(M(G)^{**})=M(G)$.\\\\

\noindent {\it{\bf Corollary 2-8.}} Suppose that $B$ is a Banach $A-bimodule$ and $B\subseteq X\subseteq B^{**}$. Let $X$ has a sequential $BAI$ $(e_n)_{n\in \mathbb{N}}\subseteq A$ and suppose that  $B$ is $WSC$. If $AX\subseteq B$, then $X=B$.\\\\

\noindent {\it{\bf Definition 2-9.}} Let $B$ be a left Banach $A-module$. Then $B^*$ is said to be left weakly completely continuous $=(\widetilde{Lwcc})$, if for each $b^\prime\in B^*$, the mapping $a\rightarrow \pi^*_\ell(b^\prime,a)$  from $A$ into $B^*$, weakly Cauchy sequence into weakly convergence ones.
If every $b^\prime\in B^*$  is  $\widetilde{Lwcc}$, then we say that $B^*$ is  $\widetilde{Lwcc}$.\\
The definition of right weakly completely continuous $(=\widetilde{Rwcc})$ is similar. We say that $b^\prime\in B^*$ is weakly completely continuous $(=\widetilde{wcc})$, if $b^\prime$ is $\widetilde{Lwcc}$ and  $\widetilde{Rwcc}$.\\\\

\noindent{\it{\bf Theorem 2-10.}} Let $B$ be a left Banach  $A-module$. Then by one of the following conditions, $B^*$ is  $\widetilde{Lwcc}$.
\begin{enumerate}
\item  $A$ is $WSC$.
\item  $B^*$ is $WSC$.
\item  $Z_{B^{**}}(A^{**})=A^{**}$. \\
\end{enumerate}
\begin{proof}
\begin{enumerate}
\item Let $(a_n)_n\subseteq A$ be weakly Cauchy sequence. Since $A$ is $WSC$, there is $a\in A$ such that $a_n\stackrel{w} {\rightarrow}a$. Now, let $b^\prime\in B^*$ and $b^{\prime\prime}\in B^{**}$. Then we have
$$\langle b^{\prime\prime},\pi^*_\ell(b^\prime,a_n)\rangle =\langle \pi^{**}_\ell(b^{\prime\prime},b^\prime),a_n\rangle \rightarrow
\langle \pi^{**}_\ell(b^{\prime\prime},b^\prime),a\rangle =\langle b^{\prime\prime},\pi^*_\ell(b^\prime,a)\rangle .$$
Thus $\pi^*_\ell(b^\prime,a_n)\stackrel{w} {\rightarrow}\pi^*_\ell(b^\prime,a)$.
\item  Proof is similar to (1).

\item Let $(a_n)_n\subseteq A$ be weakly Cauchy sequence. Since the sequence $(a_n)_n\subseteq A$ is weakly bounded in $A$, it has subsequence  such as $(a_{n_k})_k\subseteq A$  such that $a_{n_k}\stackrel{w^*} {\rightarrow}a^{\prime\prime}$ for some $a^{\prime\prime}\in A^{**}$. Then for each $b^\prime\in B^*$ and $b^{\prime\prime}\in B^{**}$, we have
$$\langle b^{\prime\prime},\pi^*_\ell(b^\prime,a_{n_k})\rangle =\langle \pi^{***}_\ell(a_{n_k},b^{\prime\prime}),b^\prime\rangle \rightarrow
\langle \pi^{***}_\ell(a^{\prime\prime},b^{\prime\prime}),b^\prime\rangle $$$$=
\langle a^{\prime\prime},\pi^{**}_\ell(b^{\prime\prime},b^\prime)\rangle =
\langle \pi^{*****}_\ell(b^{\prime\prime},b^\prime),a^{\prime\prime}\rangle $$$$=
\langle b^{\prime\prime},\pi^{****}_\ell(b^\prime,a^{\prime\prime})\rangle .$$
It is enough, we show that $\pi^{****}_\ell(b^\prime,a^{\prime\prime})\in B^*$.
Suppose that $(b^{\prime\prime}_\alpha)_\alpha \subseteq B^{**}$ such that $b^{\prime\prime}_\alpha \stackrel{w^*} {\rightarrow}b^{\prime \prime}$. Then since $Z_{B^{**}}(A^{**})=A^{**}$, we have
$$\langle \pi^{****}_\ell(b^\prime,a^{\prime\prime}),b^{\prime\prime}_\alpha\rangle
=\langle b^\prime,\pi^{***}_\ell(a^{\prime\prime},b^{\prime\prime}_\alpha)\rangle
=\langle \pi^{***}_\ell(a^{\prime\prime},b^{\prime\prime}_\alpha),b^\prime\rangle $$$$
\rightarrow \langle \pi^{***}_\ell(a^{\prime\prime},b^{\prime\prime}),b^\prime\rangle
=\langle \pi^{****}_\ell(b^\prime,a^{\prime\prime}),b^{\prime\prime}\rangle .$$
It follows that $\pi^{****}_\ell(b^\prime,a^{\prime\prime})\in (B^{**},weak^*)^*=B^*$.\\
\end{enumerate}
\end{proof}

\noindent {\it{\bf Example 2-11.}}\\ i) Suppose that $G$ is a compact group and $1\leq p\leq\infty$. Take $L^p(G)$ and $M(G)$ as $L^1(G)-bimodule$. Since $L^1(G)$ is $WSC$, by using Theorem 2-10, $L^p(G)$ and $M(G)^*$ are $\widetilde{Lwcc}$. We know that $c_0^*=\ell^1$ and $c_0$ is a $\ell^1-bimodule$. Then since $\ell^1$ is $WSC$, by using Theorem 2-10, $\ell^1$ is $\widetilde{Lwcc}$.\\
ii) Let $B$ be a reflexive Banach space. Then by using Theorem 2.6.23 from [5], $B\widehat{\otimes}B^*$, $N(B)$, the space of nucler operator on $B$, $K(B)$, the space of compact operators on $B$ and $W(B)$, the space of weakly compact operators on $B$ are Arens regular, and so by proceeding theorem, they are $\widetilde{Lwcc}$.\\\\

 Let $A$ be a Banach space. An element $a^{\prime\prime}$ of $A^{**}$ is said to be $Baire-1$ if there exists a sequence $(a_n)_n$ in $A$ that converges to $a^{\prime\prime}$ in the $weak^*$ topology of  $A^{**}$. The collection of $Baire-1$ elements of  $A^{**}$ is denoted by $\mathfrak{B}_1(A)$.\\\\

\noindent{\it{\bf Theorem 2-12.}} Let $B$ be a left Banach  $A-module$. Then $\mathfrak{B}_1(A)\subseteq Z_{B^{**}}(A^{**})$ if and only if $B^*$ is  $\widetilde{Lwcc}$.\\
\begin{proof} Let $B^*$ is  $\widetilde{Lwcc}$ and suppose that $a^{\prime\prime}\in \mathfrak{B}_1(A)$. Then there is a sequence $(a_n)_n\subseteq A$ such that $a_{n}\stackrel{w^*} {\rightarrow}a^{\prime\prime}$. It follows that
$(a_n)_n$ is weakly Cauchy sequence in $A$. Then there is $y^\prime\in B^*$ such that $\pi_\ell ^*(b^\prime ,a_{n})\stackrel{w} {\rightarrow}y^{\prime}$.\\
On the other hand, for every $b^{\prime\prime}\in B^{**}$, we have the following equality.
$$\langle b^{\prime\prime},y^\prime\rangle =\lim_n\langle b^{\prime\prime},\pi_\ell ^*(b^\prime ,a_{n})\rangle =\lim_n\langle \pi_\ell ^{**}(b^{\prime\prime},b^\prime) ,a_{n}\rangle =\lim_n\langle  a_{n},\pi_\ell ^{**}(b^{\prime\prime},b^\prime)\rangle
$$

$$=\langle  a^{\prime\prime},\pi_\ell ^{**}(b^{\prime\prime},b^\prime )\rangle =\langle  \pi_\ell ^{*****}(b^{\prime\prime},b^{\prime}),a^{\prime\prime}\rangle =\langle  b^{\prime\prime},\pi_\ell ^{****}(b^{\prime},a^{\prime\prime})\rangle .$$

It follows that $\pi_\ell ^{****}(b^{\prime},a^{\prime\prime})\in B^*$.
Suppose that $(b^{\prime\prime}_\alpha)_\alpha \subseteq B^{**}$ such that $b^{\prime\prime}_\alpha \stackrel{w^*} {\rightarrow}b^{\prime \prime}$. Then for each $b^\prime\in B^*$, we have
$$\langle \pi^{***}_\ell(a^{\prime\prime},b^{\prime\prime}_\alpha),b^\prime\rangle =
\langle \pi_\ell ^{****}(b^{\prime},a^{\prime\prime}),b^{\prime\prime}_\alpha\rangle =
\langle b^{\prime\prime}_\alpha ,\pi_\ell ^{****}(b^{\prime},a^{\prime\prime})\rangle $$$$\rightarrow
\langle b^{\prime\prime} ,\pi_\ell ^{****}(b^{\prime},a^{\prime\prime})\rangle =
\langle \pi^{***}_\ell(a^{\prime\prime},b^{\prime\prime}),b^\prime\rangle .$$
Consequently we have $a^{\prime \prime}\in Z_{B^{**}}(A^{**})$.\\
Conversely, Let $\mathfrak{B}_1(A)\subseteq Z_{B^{**}}(A^{**})$ and suppose that $(a_n)_n\subseteq A$ is weakly Cauchy sequence in $A$. Then it have subsequence such as $(a_{n_k})_k\subseteq A$  such that $a_{n_k}\stackrel{w^*} {\rightarrow}a^{\prime\prime}$ for some $a^{\prime\prime}\in A^{**}$. It follows that $a^{\prime\prime}\in \mathfrak{B}_1(A)$, and so $a^{\prime \prime}\in Z_{B^{**}}(A^{**})$. It is similar to   Theorem 2-10, we have
$\pi^{****}_\ell(b^\prime,a^{\prime\prime})\in B^*$ for some $b^\prime\in B^*$. Then for every $b^{\prime\prime}\in B^{**}$, we have
$$\lim_n\langle b^{\prime\prime},\pi_\ell^*(b^\prime,a_n)\rangle =\lim_n\langle \pi_\ell^{**}(b^{\prime\prime},b^\prime),a_n\rangle =
\lim_k\langle \pi_\ell^{**}(b^{\prime\prime},b^\prime),a_{n_k}\rangle $$$$=\lim_k\langle \pi_\ell^{***}(a_{n_k},b^{\prime\prime}),b^\prime\rangle =
\langle \pi_\ell^{***}(a^{\prime\prime},b^{\prime\prime}),b^\prime\rangle =
\langle \pi_\ell ^{*****}(b^{\prime},b^{\prime\prime}),a^{\prime\prime}\rangle $$$$=
\langle b^{\prime\prime},\pi_\ell ^{****}(b^{\prime},a^{\prime\prime})\rangle .$$
It follows that $\pi_\ell^*(b^\prime,a_n)\stackrel{w} {\rightarrow}\pi_\ell ^{****}(b^{\prime},b^{\prime\prime})$ in $B^*$.
Thus $B^*$ is  $\widetilde{Lwcc}$.\\

\end{proof}

\noindent Let $B$ be a Banach  $A-bimodule$ and $a^{\prime\prime}\in A^{**}$. We define the locally topological center of the left and right module actions of $a^{\prime\prime}$ on $B^{**}$, respectively, as follows\\
$$Z_{a^{\prime\prime}}^t(B^{**})=Z_{a^{\prime\prime}}^t(\pi_\ell^t)=\{b^{\prime\prime}\in B^{**}:~~~\pi^{t***t}_\ell(a^{\prime\prime},b^{\prime\prime})=
\pi^{***}_\ell(a^{\prime\prime},b^{\prime\prime})\},$$
$$Z_{a^{\prime\prime}}(B^{**})=Z_{a^{\prime\prime}}(\pi_r^t)=\{b^{\prime\prime}\in B^{**}:~~~\pi^{t***t}_r(b^{\prime\prime},a^{\prime\prime})=
\pi^{***}_r(b^{\prime\prime},a^{\prime\prime})\}.$$\\
It is clear that ~~~~~~~$$\bigcap_{a^{\prime\prime}\in A^{**}}Z_{a^{\prime\prime}}^t(B^{**})=Z_{A^{**}}^t(B^{**})=
Z(\pi_\ell^t),$$ ~~~~~~~ ~~\\ $$~~~~~~~~~~~~~~~~~~~~~~~~~~~~\bigcap_{a^{\prime\prime}\in A^{**}}Z_{a^{\prime\prime}}(B^{**})=Z_{A^{**}}(B^{**})=
Z(\pi_r).$$\\\\

\noindent{\it{\bf Theorem 2-13.}} Let $B$ be a left Banach  $A-module$. Then we have the following assertions.
\begin{enumerate}
\item  Suppose that $B$ has a sequential $BAI$ $(e_n)_n\subseteq A$ such that $Z_{e^{\prime\prime}}(B^{**})A\subseteq B$ where  $e^{\prime\prime}$ is a mixed unit for $A^{**}$ and $e_{n}\stackrel{w^*} {\rightarrow}e^{\prime\prime}$.  If $B$ is $WSC$, then $Z_{e^{\prime\prime}}(B^{**})=B$.
\item  If $B^*$ is $WSC$, then $\mathfrak{B}_1(A)\subseteq Z_{B^{**}}(A^{**})$.
\item  Assume that $B$ has a sequential $LBAI$ $(e_n)_n\subseteq A$ and $B^*$ is $WSC$. If $A$ is a left ideal in $A^{**}$, then $Z_{B^{**}}(A^{**})=A^{**}$.
\item  Assume that $B^*$ is $WSC$ and $A$ is a right ideal in $A^{**}$. If $Z_{e^{\prime\prime}}(A^{**})=A^{**}$, then $Z_{B^{**}}(A^{**})=A^{**}$.\\
\end{enumerate}

\begin{proof}
\begin{enumerate}
\item  Since $B\subseteq Z_{e^{\prime\prime}}(B^{**})$ for every $b\in B$, we have $\pi^{***}_r(b,e^{\prime\prime})=w^*-\lim_n \pi^{}_r(b,e_n)=b$. Let $b^{\prime\prime}\in Z_{e^{\prime\prime}}(B^{**})$. Suppose that $(b_\alpha)_\alpha \subseteq B$ such that $b_\alpha \stackrel{w^*} {\rightarrow}b^{\prime \prime}$. Then for every $b^\prime\in B^*$, we have
$$\lim_n\langle \pi^{***}_r(b^{\prime \prime},e_n),b^\prime\rangle =\langle \pi^{***}_r(b^{\prime \prime},e^{\prime\prime}),b^\prime\rangle =\lim_\alpha\langle \pi^{***}_r(b_\alpha , e^{\prime\prime}),b^\prime\rangle $$$$=\lim_{\alpha}\langle b^\prime ,b_\alpha \rangle =\langle b^{\prime \prime},b^\prime\rangle .$$
It follows that $w^*-\lim \pi^{***}_r(b^{\prime\prime},e_n)=b^{\prime\prime}$. Since $\pi^{***}_r(b^{\prime\prime},e_n)\in B$ and $B$ is $WSC$, $b^{\prime\prime}\in B$.

\item Assume that $B^*$ is $WSC$. Then by Theorem 2-10, $B^*$ is  $\widetilde{Lwcc}$, and so by using Theorem 2-11, we have $\mathfrak{B}_1(A)\subseteq Z_{B^{**}}(A^{**})$.

\item Assume that $B^*$ is $WSC$, by using part (2), we have $\mathfrak{B}_1(A)\subseteq Z_{B^{**}}(A^{**})$. Let $a^{\prime \prime}\in A^{**}$ and  suppose that  $e^{\prime \prime}\in A^{**}$ is a left unit for $A^{**}$ such that
$e_{n}\stackrel{w^*} {\rightarrow}e^{\prime\prime}$.
Then for every $a^{\prime \prime}\in A^{**}$, we have  $e_{n}a^{\prime \prime}\stackrel{w^*} {\rightarrow}e^{\prime\prime}a^{\prime \prime}=a^{\prime \prime}$. Since $AA^{**}\subseteq A$, $a^{\prime \prime} \in \mathfrak{B}_1(A)$. Consequently we have $a^{\prime \prime}\in Z_{B^{**}}(A^{**})$.
\item  Proof is similar to (3).\\
\end{enumerate}
\end{proof}

\noindent {\it{\bf Example 2-14.}} i) Let $G$ be a  compact group. We know that $L^1(G)$ has a sequential $BAI$ and  it is $WSC$ Banach algebra. Assume that  $e^{\prime\prime}$ is a mixed unit for $L^1(G)^{**}$. Since
$$Z_{e^{\prime\prime}}(L^1(G)^{**})L^1(G)\subseteq L^1(G)^{**}L^1(G)\subseteq L^1(G),$$
by using the proceeding theorem, we have
$$Z_{e^{\prime\prime}}(L^1(G)^{**})=L^1(G).$$
ii) Let $G$ be a locally compact group. In the proceeding theorem, if we take $B=c_0(G)$ and $A=\ell^1(G)$, then it is clear that $B$ is a Banach $A-bimodule$. Since $\ell^1(G)=c_0(G)^*$ is a $WSC$, $\mathfrak{B}_1(\ell^1(G))\subseteq Z^\ell_{(\ell^1(G))^{\infty}}(\ell^\infty(G)^*).$
\\\\

In the following, we extend some propositions from [21, 22] into module actions that  will be studied by Neufang.
 Let $B$ be a left Banach   $A-module$, we denote $\mathbf{L}_{A^{**}}(B^*)$ the space of linear left $A^{**}-module$ maps of $B^*$ into itself, the subspace of bounded respectively $weak^*-to-weak^*-continuous$ module maps are denote by   $\mathbf{B}_{A^{**}}(B^*)$ and $\mathbf{B}^\sigma_{A^{**}}(B^*)$, respectively. Let $B$ be a Banach space and $k\geq \it \mathcal{N}_0$ a cardinal number. A function $b^{\prime\prime}\in B^{**}$ is called $weak^*-k-continuous$ if for all nets $(b^\prime_\alpha)_{\alpha\in I}\subseteq Ball(B^*)$ of cardinality $\it \mathcal{N}_0\leq card I\leq k$ with
$b^\prime_\alpha\stackrel{w^*} {\rightarrow}0$, we have ~$\langle b^{\prime\prime},b^\prime_\alpha\rangle \rightarrow 0$. We recalled  that $B$ has the Mazur prperty of level $k$ [property $(M_k)$ ] if every $weak^*-k-continuous$ functional $b^{\prime\prime}\in B^{**}$ actually is an element of $B$, see [21, 22].\\
In the following we extend the definitions of left [left uniform] factorization property of level $k$ to module actions which introduced in [22] by Neufang for Banach algebra  and we study some problems related to module actions that they are  extension  of  Theorem 2-3 and Proposition 2-6 from [22]. The way of us prove is the same as [22].    \\\\

\noindent {\it{\bf Definition 2-15.}} Let $B$ be a right Banach  $A-module$ and $k$ be a cardinal number. We say that $B^*$ has\\
i) the left $A^{**}$ factorization property of level $k$ [property $(F_k)$] if for any family of functionals $(b^\prime_\alpha)_{\alpha\in I}\subseteq Ball(B^*)$ with $cardI=k$, there exists a family $(a^{\prime\prime}_\alpha)_{\alpha\in I}\subseteq Ball(A^{**})$ and one single functional $b^\prime\in B^*$ such that the factorization formula
$$b^\prime_\alpha=\pi_r^{**}(a^{\prime\prime}_\alpha,b^\prime)$$
hold for all $\alpha\in I$.\\
ii) the left uniform $A^{**}$ factorization property of level $k$ [property $(UF_k)$] if there exists a family $(a^{\prime\prime}_\alpha)_{\alpha\in I}\subseteq Ball(A^{**})$ with $cardI=k$ such that  for any family of functionals $(b^\prime_\alpha)_{\alpha\in I}\subseteq Ball(B^*)$ there is one single functional $b^\prime\in B^*$ such that the factorization formula
$$b^\prime_\alpha=\pi_r^{**}(a^{\prime\prime}_\alpha,b^\prime)$$
hold for all $\alpha\in I$.\\\\

\noindent {\it{\bf Theorem 2-16.}} Let $B$ be a right Banach  $A-module$ and $k$ be a cardinal number and $A$ satisfying $(M_k)$-property. Assume that $B^*$ has the property $(F_k)$. Then the following two statements hold:\\
i) $Z^\ell_{A^{**}}(B^{**})=B$.\\
ii) $\mathbf{L}_{A^{**}}(B^{*})=\mathbf{B}^\sigma_{A^{**}}(B^{*})$
\begin{proof} i) Let $b^{\prime\prime}\in Z^\ell_{A^{**}}(B^{**})$. Consider a net $(b^\prime_\alpha)_{\alpha\in I}\subseteq Ball(B^*)$, where $cardI\leq k$, which $weak^*$ converges to $0$. By property $(M_k)$, for all $a^{\prime\prime}\in A^{**}$, we only have to show that $\langle \pi^{***}_r(b^{\prime\prime},a^{\prime\prime}),b^\prime_\alpha\rangle $ converges to $0$. It suffices to prove that every convergent subnet of $(\langle \pi^{***}_r(b^{\prime\prime},a^{\prime\prime}),b^\prime_\alpha\rangle )_\alpha$ converges to $0$. Fix such a convergent subnet $(\langle \pi^{***}_r(b^{\prime\prime},a^{\prime\prime}),b^\prime_{\alpha_\beta}\rangle )_\beta$. By property $(F_k)$ for all $\alpha\in I$, we have the factorization $b^\prime_\alpha=\pi_r^{**}(a^{\prime\prime}_\alpha,b^\prime)$ where
$(a^{\prime\prime}_\alpha)_{\alpha\in I}\subseteq Ball(A^{**})$ and $b^\prime\in B^*$. Since the net
$(a^{\prime\prime}_{\alpha_\beta})_{\beta\in I}\subseteq Ball(A^{**})$ is bounded, there exists a $weak^*-convergent$ subnet $(a^{\prime\prime}_{{\alpha_\beta}_\gamma})_{\gamma\in I}\subseteq Ball(A^{**})$ such that
$ a^{\prime\prime}_{{\alpha_\beta}_\gamma}\stackrel{w^*} {\rightarrow}x^{\prime\prime}\in Ball(A^{**})$. Then for all $b\in B$, we have
$$\langle \pi^{**}_r(x^{\prime\prime},b^\prime),b\rangle =\langle x^{\prime\prime},\pi^{*}_r(b^\prime,b)\rangle
=lim_\gamma\langle a^{\prime\prime}_{{\alpha_\beta}_\gamma},\pi^{*}_r(b^\prime,b)\rangle  $$
$$=lim_\gamma\langle \pi^{**}_r(a^{\prime\prime}_{{\alpha_\beta}_\gamma},b^\prime),b)\rangle
=lim_\gamma\langle b^{\prime}_{{\alpha_\beta}_\gamma},b)\rangle =0.$$
Hence, we finally deduce , using that $b^{\prime\prime}\in Z_{A^{**}}(B^{**})$:
$$lim_\beta\langle \pi^{***}_r(b^{\prime\prime},a^{\prime\prime}),b^{\prime}_{{\alpha_\beta}}\rangle =
lim_\gamma\langle \pi^{***}_r(b^{\prime\prime},a^{\prime\prime}),b^{\prime}_{{\alpha_\beta}_\gamma}\rangle $$$$
=lim_\gamma\langle \pi^{***}_r(b^{\prime\prime},a^{\prime\prime}),\pi^{**}_r(a^{\prime\prime}_
{{\alpha_\beta}_\gamma},b^\prime)\rangle
=lim_\gamma\langle \pi^{***}_r(\pi^{***}_r(b^{\prime\prime},a^{\prime\prime}),a^{\prime\prime}_
{{\alpha_\beta}_\gamma}),b^\prime)\rangle $$
$$=\langle \pi^{***}_r(\pi^{***}_r(b^{\prime\prime},a^{\prime\prime}),x^{\prime\prime}
),b^\prime)\rangle =\langle \pi^{***}_r(b^{\prime\prime},a^{\prime\prime}),\pi^{**}_r(x^{\prime\prime},b^\prime)\rangle =0$$
which yields the desired convergence.\\
ii) Let $\phi\in \mathbf{L}_{A^{**}}{B^{*}})$ and $\phi$ is unbounded. Then, there is a sequence $(b_n^\prime)_{n\in \mathbb{N}}\subseteq Ball(B^*)$ such that $\parallel\phi(b^\prime_n)\parallel\geq n$ for all $n\in \it \mathbb{N}$.\\
Using the factorization
$$b^\prime_n=\pi^{**}_r(a^{\prime\prime}_n, b^\prime)~~~~~~(n\in  \mathbb{N})$$
where $(a_n^{\prime\prime})_n\subseteq Ball(A^{**})$ and $b^\prime\in B^*$. We obtain that for all $n\in \it N$
$$n\leq \parallel\phi(b^\prime_n)\parallel=\parallel\phi(\pi^{**}_r(a^{\prime\prime}_n, b^\prime))\parallel=\parallel \pi^{**}_r(a^{\prime\prime}_n, \phi (b^\prime))\parallel\leq \parallel\phi(b^\prime)\parallel,$$
 a contradiction.\\
 Now, we show that  $\phi\in \mathbf{B}_{A^{**}}{B^{**}})$ is automatically  $weak^*-to-weak^*-continuous$. Due to property
$(M_k)$, we only need to show that for any net $(b_\alpha^\prime)_{\alpha\in I}\subseteq Ball(B^*)$, $\it \mathcal{N}_0\leq card I\leq k$, such that $b^\prime_\alpha\stackrel{w^*} {\rightarrow}0$, we have $\phi(b^\prime_\alpha)\stackrel{w^*} {\rightarrow}0$. For all $b\in B$, we show that $\langle \phi(b^\prime_{\alpha_\beta}),b\rangle \rightarrow 0$. Property $(F_k)$ entails the factorization $b^\prime_\alpha=\pi_r^{**}(a^{\prime\prime}_\alpha,b^\prime)$
where $(a^{\prime\prime}_\alpha)_{\alpha\in I}\subseteq Ball(A^{**})$ and  $b^\prime\in B^*$. Assume that
$(a^{\prime\prime}_{{\alpha_\beta}_\gamma})_{\gamma\in I}\subseteq Ball(A^{**})$ such that
$ a^{\prime\prime}_{{\alpha_\beta}_\gamma}\stackrel{w^*} {\rightarrow}a^{\prime\prime}\in Ball(A^{**})$ and $b\in B$, then we have the following qualities
$$\langle \pi^{**}_r(a^{\prime\prime},b^\prime),b\rangle =
lim_\gamma\langle \pi^{**}_r(a^{\prime\prime}_{{\alpha_\beta}_\gamma},b^\prime),b)\rangle =
lim_\gamma\langle b^{\prime\prime}_{{\alpha_\beta}_\gamma},b)\rangle =0.$$
Hence we obtain the following results
$$lim_\beta\langle \phi(b^{\prime\prime}_{{\alpha_\beta}}),b\rangle =lim_\beta\langle \phi
(\pi^{**}_r(a^{\prime\prime}_{\alpha_\beta},b^\prime)),b\rangle =lim_\gamma
\langle \pi^{**}_r(a^{\prime\prime}_{{\alpha_\beta}_\gamma},\phi(b^\prime)),b)\rangle $$
$$=\langle \pi^{**}_r(a^{\prime\prime},\phi(b^\prime)),b\rangle =\langle \phi(\pi^{**}_r(a^{\prime\prime},b^\prime)),b\rangle =0.$$\\
\end{proof}
\noindent{\it{\bf Example 2-17.}} Let $G$ be a locally compact non-compact group with compact covering number $card(G)$. By using  [22 , Proposition 3.4],   $L^1(G)$ and $M(G)$ have Mazur property of level $cardG.\it \mathcal{N}_0$. Also by [21 , Theorem 3.2] and [22 , Theorem 3.2], respectively, $L^{\infty}(G)$ and $M(G)^*$ have $(F_k)$-property.
Since  $L^1(G)$ is a $M(G)-bimodule$ under multiplication as convolution,
$L^{\infty}(G)$ has the $M(G)^{**}$ factorization property of level $cardG.\it \mathcal{N}_0$. By using Theorem 2-16, we conclude that  $Z^\ell_{M(G)^{**}}{(L^1(G)^{**})}= {L^1(G)}$. It is similar that
$Z^\ell_{L^1(G)^{**}}{(M(G)^{**})}= {M(G)}$.\\\\

\noindent{\it{\bf Theorem 2-18.}} Suppose that $B$ is a  Banach $A-bimodule$. If
$Z^\ell_{A^{**}}(B^{**})=B$, then
 $\mathbf{B}_{A^{**}}(B^{*})=\mathbf{B}^\sigma_{A^{**}}(B^{*})$.
\begin{proof} Let $\phi\in \mathbf{B}_{A^{**}}(B^{*})$. We show that $\phi^*(B)\subseteq Z^\ell_{A^{**}}(B^{**})$. Assume that $(a^{\prime\prime}_\alpha)_\alpha\subseteq A^{**}$ such that $a^{\prime\prime}_\alpha\stackrel{w^*} {\rightarrow}0$. Then for all $b^\prime\in B^*$ and $b\in B$, we have
$$\langle \pi^{***}_r(\phi^*(b),a^{\prime\prime}_\alpha),b^\prime\rangle = \langle \phi^*(b),\pi^{**}_r(a^{\prime\prime}_\alpha,b^\prime)\rangle =
\langle b,\phi(\pi^{**}_r(a^{\prime\prime}_\alpha,b^\prime))\rangle $$
$$\langle \phi(\pi^{**}_r(a^{\prime\prime}_\alpha,b^\prime)),b\rangle =\langle \pi^{**}_r(a^{\prime\prime}_\alpha,\phi(b^\prime)),b\rangle
=\langle a^{\prime\prime}_\alpha,\pi^{*}_r(\phi(b^\prime),b)\rangle \rightarrow 0.$$
We conclude that $\phi^*(b)\in Z^\ell_{A^{**}}(B^{**})=B$ and result is hold.\\
\end{proof}
\noindent{\it{\bf Example 2-19.}} Suppose that $G$ is either an infinite locally compact group such that $k(G)\geq 2^{b(G)}$
where $k(G)$ is the compact covering number, that is, the least cardinality of a covering of $G$ by compact subsets and $b(G)$ is the local weight of $G$, that is, the least cardinality of an open basis at the neutral element of $G$ or  $G$ is locally compact non-compact group with non-measurable cardinality. By using  [21, Proposition 3.5], [22, Proposition 3.4] and Theorem 2-18, we conclude that
$$\mathbf{L}_{L^1(G)^{**}}({M(G)^{**}})=\mathbf{B}^\sigma_{L^1(G)^{**}}({M(G)^{**}}).$$\\

\begin{center}
\section{ \bf Cohomological groups of Banach algebras  }
\end{center}
Let $B$ be a Banach $A-bimodule$ and let $n\geq 1$. Suppose that $B^{(n)}$ is an  $n-th~dual$ of  $B$. Then $B^{(n)}$ is also Banach $A-bimodule$, that is, for every $a\in A$, $b^{(n)}\in B^{(n)}$ and $b^{(n-1)}\in B^{(n-1)}$, we define
 $$\langle b^{(n)}a,b^{(n-1)}\rangle= \langle b^{(n)},ab^{(n-1)}\rangle,$$
 and
 $$\langle ab^{(n)},b^{(n-1)}\rangle= \langle b^{(n)},b^{(n-1)}a\rangle.$$
 In the following theorem, we extend  Theorem 1.9 from [6] into general situation and  we give similar proof.\\\\

\noindent{\it{\bf Theorem 3-1.}} Let $B$ be a Banach $A-bimodule$ and let $n\geq 1$. If $H^1(A,B^{(n+2)})=0$, then $H^1(A,B^{(n)})=0$.\\
\begin{proof} Let $D\in Z^1(A,B^{(n)})$ and suppose that  $i: B^{(n)}\rightarrow B^{(n+2)}$ is the canonical linear mapping as $A-bimodule$ homomorphism. Take $\widetilde{D}=ioD$. Then we can be viewed $\widetilde{D}$ as an element of
$Z^1(A,B^{(n+2)})$. Since $H^1(A,B^{(n+2)})=0$, there exist $b^{(n+2)}\in B^{(n+2)}$ such that for every $a\in A$, we have $$\widetilde{D}(a)=ab^{(n+2)}-b^{(n+2)}a.$$
Set a $A-linear$ mapping $P$ from $B^{(n+2)}$ into $B^{(n)}$ such that $Poi=I_{B^{(n)}}$. Then we have $Po\widetilde{D}=(Poi)oD=D$, and so for every $a\in A$, we conclude that $D(a)=Po\widetilde{D}(a)=aP(b^{(n+2)})-P(b^{(n+2)})a$. It follows that $D\in N^1(A,B^{(n)})$. Consequently we have $H^1(A,B^{(n)})=0$.\\
\end{proof}

In every parts of this section, for left or right Banach $A-module$ $B$, we take $\pi_\ell(a,b)=ab$ and  $\pi_r(b,a)=ba$, respectively,  for each $a\in A$ and $b\in B$.\\
Let $A^{(n)}$ and  $B^{(n)}$  be $n-th~dual$ of $A$ and $B$, respectively. By [28, page 4132-4134], if $n\geq 0$ is an even number, then  $B^{(n)}$ is a Banach $A^{(n)}-bimodule$. Then for $n\geq 2$,   we define  $B^{(n)}B^{(n-1)}$ as a subspace of $A^{(n-1)}$, that is, for all $b^{(n)}\in B^{(n)}$,  $b^{(n-1)}\in B^{(n-1)}$ and  $a^{(n-2)}\in A^{(n-2)}$ we define
$$\langle b^{(n)}b^{(n-1)},a^{(n-2)}\rangle=\langle b^{(n)},b^{(n-1)}a^{(n-2)}\rangle.$$
If $n$ is odd number, then for $n\geq 1$,   we define  $B^{(n)}B^{(n-1)}$ as a subspace of $A^{(n)}$, that is, for all $b^{(n)}\in B^{(n)}$,  $b^{(n-1)}\in B^{(n-1)}$ and  $a^{(n-1)}\in A^{(n-1)}$ we define
$$\langle b^{(n)}b^{(n-1)},a^{(n-1)}\rangle=\langle b^{(n)},b^{(n-1)}a^{(n-1)}\rangle.$$
and if $n=0$, we take $A^{(0)}=A$ and $B^{(0)}=B$.\\
We also define the topological centers of module actions of $A^{(n)}$ on  $B^{(n)}$ as follows

$${Z}^\ell_{A^{(n)}}(B^{(n)})=\{b^{(n)}\in B^{(n)}:~the~map~~a^{(n)}\rightarrow b^{(n)} a^{(n)}~:~A^{(n)}\rightarrow B^{(n)}$$$$~is~~~weak^*-to-weak^*~continuous\}$$
$${Z}^\ell_{B^{(n)}}(A^{(n)})=\{a^{(n)}\in A^{(n)}:~the~map~~b^{(n)}\rightarrow a^{(n)} b^{(n)}~:~B^{(n)}\rightarrow B^{(n)}$$$$~is~~~weak^*-to-weak^*~continuous\}.$$\\

In the following, for every $n\geq 1$, by using  Arens regularity of left module action of $A^{(2n)}$ on $B^{(2n)}$, for every derivation $D:A\rightarrow B^{(2n)}$, we find a derivation $\widetilde{D}$ from $A^{(2n)}$ into $B^{(2n)}$ such that $\widetilde{D}(a)=D(a)$ for every $a\in A$. We have some conclusions in cohomological groups.\\\\

\noindent{\it{\bf Theorem 3-2.}} Let $B$ be a Banach $A-bimodule$ and  $D:A\rightarrow B^{(2n)}$ be a continuous derivation. Assume that $Z^\ell_{A^{(2n)}}(B^{(2n)})=B^{(2n)}$. Then there is a continuous derivation $\widetilde{D}: A^{(2n)}\rightarrow B^{(2n)}$ such that $\widetilde{D}(a)=D(a)$ for every $a\in A$.
\begin{proof} By using [6, Proposition 1.7], the linear mapping $D^{\prime\prime}:A^{**}\rightarrow B^{(2n+2)}$ is a continuous derivation. Take $X=B^{(2n-2)}$. Since $Z_{A^{(2n)}}(X^{**})=Z_{A^{(2n)}}(B^{(2n)})=B^{(2n)}=X^{**}$, by [6, Proposition 1.8], the canonical projection $P:  X^{(4)}\rightarrow X^{**}$ is a $A^{**}-bimodule $ morphism. Set $\widetilde{D}=PoD^{\prime\prime}$. Then $\widetilde{D}$ is a continuous derivation from $A^{**}$ into $B^{(2n)}$. Now by replacing $A^{**}$ by $A$ and repeating of the proof, we obtain the result. \\ \end{proof}

\noindent{\it{\bf Corollary 3-3.}} Let $B$ be a Banach $A-bimodule$ and $n\geq 0$. If $Z^\ell_{A^{(2n)}}(B^{(2n)})=B^{(2n)}$ and $H^1(A^{(2n+2)},B^{(2n+2)})=0$, then  $H^1(A,B^{(2n)})=0$.
\begin{proof} By using [6, Proposition 1.7] and proceeding theorem the result is hold.\\ \end{proof}

Let $A$ be a Banach algebra and $n\geq 0$. Then $A$ is $n-weakly$ amenable if $H^1(A,A^{(n)})=0$. $A$ is permanently weakly amenable if $A$ is $n-weakly$ amenable for each $n\geq 0$.\\\\

\noindent{\it{\bf Corollary 3-4}} [6]. Let $A$ be a Banach algebra such that $A^{(2n)}$ is Arens regular and
 $H^1(A^{(2n+2)}),A^{(2n+2)})=0$ for each $n\geq 0$. Then $A$ is $2n-weakly$ amenable.\\\\

Assume that $A$ is Banach algebra and $n\geq 0$. We define $A^{[n]}$ as a subset of $A$ as follows
$$A^{[n]}=\{a_1a_2...a_n:~a_1,a_2,...a_n\in A\}.$$
We write $A^n$  the linear span of $A^{[n]}$ as a subalgebra of $A$.\\\\

\noindent{\it{\bf Theorem 3-5.}} Let $A$ be a Banach algebra and $n\geq 0$. Let $A^{[2n]}$ dense in $A$ and suppose that $B$ is a Banach $A-bimodule$. Assume that $AB^{**}$ and $B^{**}A$ are subsets of $B$. If $H^1(A,B^*)=0$, then  $H^1(A ,B^{(2n+1)})=0$.
\begin{proof} For $n=0$ the result is clear.  Let $B^\perp$ be the space of functionals in $A^{(2n+1)}$ which annihilate $i(B)$ where $i:B\rightarrow
B^{(2n)}$ is a natural canonical mapping. Then, by using [29, lemma 1], we have the following equality.
$$B^{(2n+1)}=i(B^*)\oplus B^\bot ,$$
as Banach $A-bimodules$, and so
$$H^1(A,B^{(2n+1)})=H^1(A,i(B^*))\oplus H^1(A,B^\bot ) .$$ Without lose generality, we replace $i(B^*)$ by $B^*$.
Since $H^1(A,B^*)=0$, it is enough to show that $H^1(A,B^\bot )=0$.\\
Now, take the linear mappings $L_a$ and $R_a$ from $B$ into itself by $L_a(b)=ab$ and $R_a(b)=ba$ for every $a\in A$.
Since $AB^{**}\subseteq B$ and $B^{**}A\subseteq B$,  $L^{**}_a(b^{\prime\prime})=ab^{\prime\prime}$ and $R^{**}_a(b^{\prime\prime})=b^{\prime\prime}a$ for every $a\in A$, respectively. Consequently, $L_a$ and $R_a$ from $B$ into itself are weakly compact. It follows that for each $a\in A$ the linear mappings $L^{(2n)}_a$ and   $R^{(2n)}_a$ from $B^{(n)}$ into $B^{(n)}$ are weakly compact and for every $b^{(2n)}\in B^{(2n)}$, we have $L^{(2n)}_a(b^{(2n)})=ab^{(2n)}\in B^{(2n-2)}$ and $R^{(2n)}_a(b^{(2n)})=b^{(2n)}a\in B^{(2n-2)}$. Set $a_1, a_2, ...,a_n\in A$ and
$b^{(2n)}\in B^{(2n)}$. Then $a_1 a_2 ...a_nb^{(2n)}$ and $b^{(2n)}a_1 a_2 ...a_n$ are belong to $B$.
Suppose that $D\in Z^1(A,B^\bot )$ and let $a, x\in A^{[n]}$. Then for every $b^{(2n)}\in B^{(2n)}$, since $xb^{(2n)}, b^{(2n)}a \in B$, we have the following equality
$$\langle D(ax), b^{(2n)}\rangle=\langle aD(x), b^{(2n)}\rangle+\langle D(a)x, b^{(2n)}\rangle$$$$=\langle D(x), b^{(2n)}a\rangle+\langle D(a), xb^{(2n)}\rangle=0.$$
It follows that $D\mid_{A^{[2n]}}$. Since $A^{[2n]}$ dense in $A$, $D=0$. Hence $H^1(A,B^\bot )=0$ and result follows.\\
\end{proof}

\noindent{\it{\bf Corollary 3-6.}}\\ i) Let $A$ be a Banach algebra with bounded left approximate identity [=$LBAI$], and let  $B$ be a Banach $A-bimodule$.  Suppose that $AB^{**}$ and $B^{**}A$ are subset of $B$. Then if $H^1(A,B^*)=0$, it follows that  $H^1(A,B^{(2n+1)})=0$.\\
ii) Let $A$ be an amenable Banach algebra and $B$ be a Banach $A-bimodule$. If $AB^{**}$ and $B^{**}A$ are subset of $B$, then $H^1(A,B^{(2n+1)})=0$.\\\\

\noindent{\it{\bf Example 3-7.}} Assume that $G$ is a compact group. Then \\
i) we know that $L^1(G)$ is $M(G)-bimodule$  and $L^1(G)$ is an ideal in the second dual of $M(G)$, $M(G)^{**}$. By using [19, corollary 1.2], we have $H^1(L^1(G),M(G))=0$. Then  for every $n\geq 1$, by using proceeding corollary, we conclude that
$$H^1(L^1(G),M(G)^{(2n+1)})=0.$$
ii) we have $L^1(G)$ is an ideal in its second dual , $L^1(G)^{**}$. By using [16], we know that $L^1(G)$ is a weakly amenable. Then by proceeding corollary, $L^1(G)$ is $(2n+1)-weakly$  amenable.\\\\

\noindent{\it{\bf Corollary 3-8.}} Let $A$ be a Banach algebra and let $A^{[2n]}$ be dense in $A$. Suppose that $AB^{**}$ and $B^{**}A$ are subset of $B$. Then the following are equivalent.
\begin{enumerate}
\item ~$H^1(A,B^*)=0$.
\item ~$H^1(A,B^{(2n+1)})=0$ for some $n\geq 0$.

\item ~~$H^1(A,B^{(2n+1)})=0$ for each $n\geq 0$.\\\\
\end{enumerate}

\noindent{\it{\bf Corollary 3-9}} [6]. Let $A$ be a weakly amenable Banach algebra such that $A$ is an ideal in $A^{**}$. Then $A$ is $(2n+1)-weakly$ amenable for each $n\geq 0$.
\begin{proof} By using Proposition 1.3 from [6] and proceeding theorem, result is hold.\\ \end{proof}

Let $B$ be a dual Banach algebra and $dimB< \infty$. Then by using Proposition 2.6.24 from [5], we know that $\mathcal{N}(B)$, the collection of all operators from $B$ into $B$, is an ideal in $\mathcal{N}(B)^{**}$. By using Corollary 4.3.6 from [24], $N(B)$ is amenable, and so it also is weakly amenable. Consequently, by using the proceeding corollary, $\mathcal{N}(B)$ is $(2n+1)-weakly$ amenable for every $n\geq 0$.\\
We know that every von Neumann algebra $A$ is an weakly amenable Banach algebra, see [24]. Now if $A$ is an ideal in its second dual, $A^{**}$, then by using proceeding corollary, $A$ is $(2n+1)-weakly$ amenable Banach algebra for each $n\geq 0$.\\

Assume that $A$ and $B$ are Banach algebra. Then $A\oplus B$ , with norm $$\parallel(a,b)\parallel=\parallel a\parallel+\parallel b\parallel,$$ and product $(a_1,b_1)(a_2,b_2)=(a_1a_2, b_1b_2)$ is a Banach algebra. It is clear that if $X$ is a Banach $A~and ~B-bimodule$, then $X$ is a Banach $A\oplus B-bimodule$.\\
In the following, we investigated the relationships between the cohomological group of  $A\oplus B$ and cohomological groups of $A$ and $B$.\\\\

\noindent{\it{\bf Theorem 3-10.}} Suppose that $A$ and $B$ are Banach algebras. Let $X$ be a Banach $A~and ~B-bimodule$. Then,   $H^1(A\oplus B, X)=0$ if and only if $H^1(A,X)=H^1(B,X)=0$.
\begin{proof} Suppose that $H^1(A\oplus B, X)=0$. Assume that $D_1\in Z^1(A,X)$ and $D_2\in Z^1(B,X)$. Take $D=(D_1,D_2)$. Then for every $a_1,a_2\in A$ and $b_1, b_2\in B$, we have
$$D((a_1,b_1)(a_2,b_2))=D(a_1a_2,b_1b_2)=(D_1(a_1a_2),D_2(b_1b_2))$$
$$=(a_1D_1(a_2)+D_1(a_1)a_2,b_1D_2(b_2)+D_2(b_1)b_2)$$
$$=(a_1D_1(a_2),b_1D_2(b_2))+(D_1(a_1)a_2+D_2(b_1)b_2)$$
$$=(a_1,b_1)(D_1(a_2),D_2(b_2))+(D_1(a_1),D_2(b_1))(a_2,b_2)$$
$$=(a_1,b_1)D(a_2,b_2)+D(a_1,b_1)(a_2,b_2).$$
It follows that $D\in Z^1(A\oplus B,X)$.
Since $H^1(A\oplus B, X)=0$, there is $x\in X$ such that $D=\delta_{x}$  where  $\delta_{x}\in N^1(A\oplus B,X)$. Since
$\delta_{x}=(\delta^1_{x},\delta^2_{x})$ where $\delta^1_{x}\in N^1(A,X)$ and $\delta^2_{x}\in N^1(B,X)$,  we have $D_1=\delta^1_{x}$ and $D_2=\delta^2_{x}$. Thus we have
$H^1(A,X)=H^1(B,X)=0$.\\
For the converse, take $A$ as an ideal in $A\oplus B$, and so by using Proposition 2.8.66 from [5], proof hold.\\ \end{proof}

Let $G$ be a locally compact group and $X$ be a Banach $L^1(G)-bimodule$. Then by [6, pp.27 and 28], $X^{**}$ is a Banach $L^1(G)^{**}-bimodule$.  Since $L^1(G)^{**}=LUC(G)^*\oplus LUC(G)^\bot$, by using proceeding theorem, we have  $H^1(L^1(G)^{**},X^{**})=0$ if and only if $H^1(LUC(G)^*,X^{**})=H^1(LUC(G)^\bot,X^{**})=0$.\\
On the other hand, we know that $L^1(G)^{**}=L^1(G)\oplus C_0(G)^\perp$. By [16], we know that, $H^1(L^1(G), L^\infty(G))=0$. By using proceeding theorem,  $$H^1(L^1(G)^{**}, L^\infty(G))=0,$$ if and only if  $H^1(C_0(G)^\perp, L^\infty(G))=0$.\\\\

\noindent{\it{\bf Corollary 3-11.}} \\
i) Suppose that $A$ and $B$ are Banach algebras. Let $X$ be a Banach $A~and ~B-bimodule$. Then   $A\oplus B$ is an amenable Banach algebra if and only if $A$ and $B$ are amenable Banach algebras.\\
ii) Let $A$ be a Banach algebra and $n\geq 1$. Then  $H^1(\oplus_{i=1}^nA,A^*)=0$ if and only if $A$ is weakly amenable.\\\\

Assume that $B$ is a Banach $A-bimodule$. In the following, we will study the relationships between two cohomological groups $H^1(A,B^*)$ and  $H^1(A^{**},B^{***})$. If  $H^1(A,B^*)=0$, we want to know that dose $H^1(A^{**},B^{***})=0$ and when its converse hold? If we give a positive answer to these questions, then we want to establish the results of them in the weak amenability of Banach algebras. As some applications of these discussion in the group algebras, for a compact group $G$, we show that every $weak^*-to-weak$ continuous derivation from $L^1(G)^{**}$ into $M(G)$ is inner. \\\\

\noindent{\it{\bf Theorem 3-12.}} Let $B$ be a Banach $A-bimodule$ and suppose that every derivation from $A$ into $B^*$ is weakly compact. If $Z^\ell_{A^{**}}(B^{**})=B^{**}$ and $H^1(A^{**},B^{***})=0$, then $H^1(A,B^*)=0$.
\begin{proof} Suppose that $D\in Z^1(A,B^*)$. Since $D:A\rightarrow B^*$ is weakly compact $D^{\prime\prime}(A^{**})\subseteq B^*$. Assume that $a^{\prime\prime}, x^{\prime\prime}\in A^{**}$ and $(a_\alpha)_\alpha, (x_\beta)_\beta\subseteq A$ such that $a_\alpha^{} \stackrel{w^*} {\rightarrow} a^{\prime\prime}$ and $x_\beta \stackrel{w^*} {\rightarrow} x^{\prime\prime}$ in $A^{**}$. Then, since $Z^\ell_{A^{**}}(B^{**})=B^{**}$, for every $b^{\prime\prime}\in B^{**}$, we have
$$\lim_\alpha \langle a_\alpha D^{\prime\prime}(x^{\prime\prime}),b^{\prime\prime}\rangle=
\lim_\alpha \langle  D^{\prime\prime}(x^{\prime\prime}),b^{\prime\prime}a_\alpha\rangle=
\langle b^{\prime\prime}a^{\prime\prime}, D^{\prime\prime}(x^{\prime\prime})\rangle$$$$=
 \langle b^{\prime\prime},a^{\prime\prime} D^{\prime\prime}(x^{\prime\prime})\rangle=
 \langle a^{\prime\prime} D^{\prime\prime}(x^{\prime\prime}),b^{\prime\prime}\rangle.$$
Consequently, we have
$$D^{\prime\prime}(a^{\prime\prime}x^{\prime\prime})=\lim_\alpha\lim_\beta D(a_\alpha x_\beta)=
\lim_\alpha\lim_\beta (D(a_\alpha) x_\beta +a_\alpha D(x_\beta))$$$$=D^{\prime\prime}(a^{\prime\prime})x^{\prime\prime}
+a^{\prime\prime}D^{\prime\prime}(x^{\prime\prime}).$$
It follows that $D^{\prime\prime}\in Z^1(A^{**},B^{***})$. Since $H^1(A^{**},B^{***})=0$, there is $b^{\prime\prime\prime}\in  B^{***}$ such that $D^{\prime\prime}=\delta_{b^{\prime\prime\prime}}$. Then we have $D=D^{\prime\prime}\mid_A=\delta_{b^{\prime\prime\prime}}\mid_A$. So for every $a\in A$, we have $D(a)=\delta_{b^{\prime\prime\prime}}(a)=ab^{\prime\prime\prime}-b^{\prime\prime\prime}a$. If we take $b^\prime=b^{\prime\prime\prime}\mid_{B}$, then $D(a)=ab^\prime-b^\prime a=\delta_{b^\prime}(a)$. We conclude that $H^1(A,B^*)=0$.\\ \end{proof}

Let $A$ be a Banach algebra, and let $B$ be a Banach $A-bimodule$. A continuous derivation $D:A\rightarrow B$ is approximately inner [resp. weakly approximately inner], if there exists a bounded net $(b_\alpha)_\alpha \subseteq B$ such that $D(a)=\lim_\alpha (ab_\alpha-b_\alpha a)$  in $B$ [resp. $D(a)=w-\lim_\alpha (ab_\alpha-b_\alpha a)$ in $B$].\\
Assume that $(b^\prime_\alpha)_\alpha \subseteq B^*$ such that $b^\prime_\alpha \stackrel{w^*} {\rightarrow} b^{\prime\prime\prime}$ in $B^{***}$. Then by assumptions of proceeding theorem, for every $a\in A$ and for every derivation from $A$ into $B^*$,   we have
$D(a)=D^{\prime\prime}\mid_A(a)=\delta_{b^{\prime\prime\prime}}\mid_A(a)=ab^{\prime\prime\prime}-b^{\prime\prime\prime
}a=w^*-\lim_\alpha (ab^\prime_\alpha-b^\prime_\alpha a)$ in $B^{***}$. Since $D(a)\in B^*$, $D(a)=w-\lim_\alpha (ab^\prime_\alpha-b^\prime_\alpha a)$ in $B^*$, and so $D$ is weakly approximately inner in $B^*$.\\\\

\noindent{\it{\bf Theorem 3-13.}}  Let $B$ be a Banach $A-bimodule$ and suppose that $D^{\prime\prime}(A^{**})\subseteq B^*$. If $Z^\ell_{A^{**}}(B^{**})=B^{**}$ and $H^1(A^{**},B^{***})=0$, then $H^1(A,B^*)=0$.
\begin{proof} Proof is similar to proceeding theorem.\\ \end{proof}

A functional $a^\prime$ in $A^*$ is said to be $wap$ (weakly almost
 periodic) on $A$ if the mapping $a\rightarrow a^\prime a$ from $A$ into
 $A^{*}$ is weakly compact.  In  [23], Pym showed that  this definition to the equivalent following condition\\
 For any two net $(a_{\alpha})_{\alpha}$ and $(b_{\beta})_{\beta}$
 in $\{a\in A:~\parallel a\parallel\leq 1\}$, we have\\
$$\\lim_{\alpha}\\lim_{\beta}\langle a^\prime,a_{\alpha}b_{\beta}\rangle=\\lim_{\beta}\\lim_{\alpha}\langle a^\prime,a_{\alpha}b_{\beta}\rangle,$$
whenever both iterated limits exist. The collection of all $wap$
functionals on $A$ is denoted by $wap(A)$. Also we have
$a^{\prime}\in wap(A)$ if and only if $\langle a^{\prime\prime}b^{\prime\prime},a^\prime\rangle=\langle a^{\prime\prime}ob^{\prime\prime},a^\prime\rangle$ for every $a^{\prime\prime},~b^{\prime\prime} \in
A^{**}$. \\
Let $B$ be a Banach left $A-module$. Then, $b^\prime\in B^*$ is said to be left weakly almost periodic functional if the set $\{\pi^*_\ell(b^\prime,a):~a\in A,~\parallel a\parallel\leq 1\}$ is relatively weakly compact. We denote by $wap_\ell(B)$ the closed subspace of $B^*$ consisting of all the left weakly almost periodic functionals in $B^*$.\\
The definition of the right weakly almost periodic functional ($=wap_r(B)$) is the same.\\
By  [23], $b^\prime\in wap_\ell(B)$ is equivalent to the following $$\langle \pi_\ell^{***}(a^{\prime\prime},b^{\prime\prime}),b^\prime\rangle=
\langle \pi_\ell^{t***t}(a^{\prime\prime},b^{\prime\prime}),b^\prime\rangle$$
for all $a^{\prime\prime}\in A^{**}$ and $b^{\prime\prime}\in B^{**}$.
Thus, we can write \\
$$wap_\ell(B)=\{ b^\prime\in B^*:~\langle \pi_\ell^{***}(a^{\prime\prime},b^{\prime\prime}),b^\prime\rangle=
\langle \pi_\ell^{t***t}(a^{\prime\prime},b^{\prime\prime}),b^\prime\rangle~~$$$$for~~all~~a^{\prime\prime}\in A^{**},~b^{\prime\prime}\in B^{**}\}.$$\\\\

\noindent{\it{\bf Corollary 3-14.}}  Let $B$ be a Banach $A-bimodule$ and let every derivation\\ $D:A\rightarrow B^*$,
satisfies $D^{\prime\prime}(A^{**})\subseteq wap_\ell(B)$. If $H^1(A^{**},B^{***})=0$, then\\ $H^1(A,B^*)=0$.\\

In Corollary 3-14, if we take $B=A$, then we obtain Theorem 2.1 from [10].\\

\noindent{\it{\bf Corollary 3-15.}}  Let $B$ be a Banach $A-bimodule$ and suppose that for every derivation $D:A\rightarrow B^*$, we have $D^{\prime\prime}(A^{**})B^{**}\subseteq A^*$. If $H^1(A^{**},B^{***})=0$, then $$H^1(A,B^*)=0.$$\\

\noindent{\it{\bf Theorem 3-16.}}  Let $B$ be a Banach $A-bimodule$ and suppose that $AA^{**}\subseteq A$ and
$B^{**}A=B^{**}$. If $H^1(A^{**},B^{***})=0$, then $H^1(A,B^*)=0$.
\begin{proof} Let   $(a_{\alpha}^{\prime\prime})_{\alpha}\subseteq A^{**}$ such that  $a^{\prime\prime}_{\alpha} \stackrel{w^*} {\rightarrow}a^{\prime\prime}$ in $A^{**}$. Assume that $b^{\prime\prime}\in B^{**}$. Since $B^{**}A=B^{**}$, there are $a\in A$ and $y^{\prime\prime}\in B^{**}$ such that $b^{\prime\prime}=y^{\prime\prime}a$. We know that $a\in A\subseteq Z^\ell_{B^{**}}(A^{**})$, and so $aa^{\prime\prime}_{\alpha} \stackrel{w^*} {\rightarrow}aa^{\prime\prime}$ in $A^{**}$. Since $(aa^{\prime\prime}_{\alpha})_{\alpha}\subseteq A$ and $aa^{\prime\prime}\in A$, we have $aa^{\prime\prime}_{\alpha} \stackrel{w} {\rightarrow}aa^{\prime\prime}$ in $A$.
Then for every $a^{\prime\prime}\in A^{**}$, we have the following equalities
$$\langle D^{\prime\prime}(x^{\prime\prime})b^{\prime\prime},a^{\prime\prime}_{\alpha}\rangle=
\langle D^{\prime\prime}(x^{\prime\prime})y^{\prime\prime}a,a^{\prime\prime}_{\alpha}\rangle=
\langle D^{\prime\prime}(x^{\prime\prime})y^{\prime\prime},aa^{\prime\prime}_{\alpha}\rangle\rightarrow
\langle D^{\prime\prime}(x^{\prime\prime})y^{\prime\prime},aa^{\prime\prime}\rangle$$$$=
\langle D^{\prime\prime}(x^{\prime\prime})b^{\prime\prime},a^{\prime\prime}\rangle.$$
It follows that $D^{\prime\prime}(x^{\prime\prime})b^{\prime\prime}\in (A^{**},weak^*)^*=A^*$. So, by using Corollary 3-15, proof  hold.\\ \end{proof}

\begin{center}
\section{ \bf The $weak^*-to- weak^*$ continuity  of derivations and Connes-amenability of Banach algebras  }
\end{center}

\noindent{\it{\bf Theorem 4-1.}} Let $B$ be a Banach $A-bimodule$  and let every derivation $D:A^{**}\rightarrow B^*$ is $weak^*-to-weak^*$ continuous. Then if $Z^\ell_{B^{**}}(A^{**})=A^{**}$ and $H^1(A,B^*)=0$, it follows that  $H^1(A^{**},B^*)=0$.
\begin{proof} Let $D:A^{**}\rightarrow B^*$ be a derivation. Then $D\mid_A:A\rightarrow B^*$ is a derivation. Since $H^1(A,B^*)=0$, there is $b^\prime\in B^*$ such that $D\mid_A=\delta_{b^\prime}$. Suppose that $a^{\prime\prime}\in A^{**}$ and
$(a_\alpha^{})_\alpha\subseteq A^{}$ such that $a_\alpha^{} \stackrel{w^*} {\rightarrow} a^{\prime\prime}$ in $A^{**}$. Then
$$D(a^{\prime\prime})=w^*-\lim_\alpha D\mid_A(a_\alpha)=w^*-\lim_\alpha\delta_{b^\prime}(a_\alpha)=w^*-\lim_\alpha(a_\alpha {b^\prime}-{b^\prime}a_\alpha)$$$$=a^{\prime\prime}{b^\prime}-{b^\prime}a^{\prime\prime}.$$
Now, we show that $b^\prime a^{\prime\prime}\in B^*$. Assume that  $(b^{\prime\prime}_{\beta})_{\beta}\in B^{**}$ such that $b^{\prime\prime}=w^*-\lim_{\beta}b^{\prime\prime}_{\beta}$. Then since $Z^\ell_{B^{**}}(A^{**})=A^{**}$, we have
$$\langle b^\prime a^{\prime\prime},b^{\prime\prime}_{\beta}\rangle=\langle  a^{\prime\prime}b^{\prime\prime}_{\beta},b^\prime\rangle\rightarrow
\langle  a^{\prime\prime}b^{\prime\prime},b^\prime\rangle=\langle b^\prime a^{\prime\prime},b^{\prime\prime}\rangle.$$
Thus, we conclude that $b^\prime a^{\prime\prime}\in (B^{**},weak^*)^*=B^*$, and so
$H^1(A^{**},B^*)=0$.\\ \end{proof}

\noindent{\it{\bf Corollary 4-2.}} Let $A$ be a Arens regular  Banach algebra and let every derivation $D:A^{**}\rightarrow A^*$ is $weak^*-to-weak^*$ continuous. Then if $A$ is weakly amenable, it follows that  $H^1(A^{**},A^*)=0$.\\\\

Let $B$  be a dual Banach space, with predual  $X$ and suppose that
$$X^\perp=\{x^{\prime\prime\prime}:~x^{\prime\prime\prime}\mid_X=0~~where~~x^{\prime\prime\prime}\in X^{***}\}=\{b^{\prime\prime}:~b^{\prime\prime}\mid_X=0~~where~~b^{\prime\prime}\in B^{**}\}.$$
Then the canonical projection $P:X^{***}\rightarrow X^*$ gives a continuous linear map $P:B^{**}\rightarrow B$. Thus, we can write the following equality
$$B^{**}=X^{***}=X^*\oplus ker P=B\oplus X^\perp ,$$
as a direct sum of Banach $A-bimodule$.\\\\

\noindent{\it{\bf Theorem 4-3.}}  Let $B$ be a Banach $A-bimodule$ and  every derivation $A^{**}$ into  $ B$ is $weak^*-to-weak$ continuous. Let $A^{**}B,BA^{**}\subseteq B$. Then  the following assertions are hold.
\begin{enumerate}
\item  If $H^1(A,B)=0$, then $H^1(A^{**},B)=0$.

\item  Suppose that $A$ has a $LBAI$. Let $B$ has a predual $X$ and let $AB^*,~B^*A\subseteq X$. Then, if $H^1(A,B)=0$, it follows that  $H^1(A^{**},B^{**})=0$.

\end{enumerate}
\begin{proof}
\begin{enumerate}
\item  Proof is similar to Theorem 4-1.

\item  Set $B^{**}=B\oplus X^\perp .$ Then we have
$$H^1(A^{**},B^{**})=H^1(A^{**},B)\oplus H^1(A^{**},X^\perp ).$$
Since $H^1(A,B)=0$, by part (1), $H^1(A^{**},B)=0$. Now let $\widetilde{D}\in Z^1(A^{**},X^\perp)$ and we take $D=\widetilde{D}\mid_A$. It is clear that $D\in Z^1(A^{**},X^\perp)$. Assume that $a^{\prime\prime}, x^{\prime\prime}\in A^{**}$ and $(a_\alpha)_\alpha, (x_\beta)_\beta\subseteq A$ such that $a_\alpha^{} \stackrel{w^*} {\rightarrow} a^{\prime\prime}$ and $x_\beta \stackrel{w^*} {\rightarrow} x^{\prime\prime}$ on $A^{**}$. Since $AB^*,~B^*A\subseteq X$, for every $b^ \prime \in B^*$, by using the $weak^*-to-weak$ continuity  of $\widetilde{D}$, we have
$$\langle \widetilde{D}(a^{\prime\prime}x^{\prime\prime}), b^\prime\rangle=\lim_\alpha\lim_\beta \langle D(a_\alpha x_\beta), b^\prime\rangle$$$$=
\lim_\alpha\lim_\beta \langle (D(a_\alpha) x_\beta +a_\alpha D(x_\beta)),b^\prime\rangle$$$$=
\lim_\alpha\lim_\beta \langle D(a_\alpha) x_\beta ,b^\prime\rangle+\lim_\alpha\lim_\beta \langle a_\alpha D(x_\beta),b^\prime\rangle$$
$$=
\lim_\alpha\lim_\beta \langle D(a_\alpha), x_\beta b^\prime\rangle+\lim_\alpha\lim_\beta \langle D(x_\beta)),b^\prime a_\alpha \rangle$$

$$=0.$$
Since $A$ has a $LBAI$, $A^{**}$ has a left unit $e^{\prime\prime}$ with respect to the first Arens product. Then  $D(x^{\prime\prime})= D(e^{\prime\prime}x^{\prime\prime})=0$, and so $D=0$.\\ \end{enumerate}
\end{proof}

\noindent{\it{\bf Example 4-4.}}\\
 i) Assume that $G$ is a compact group. Then we know that $L^1(G)$ is $M(G)-bimodule$  and $L^1(G)$ is an ideal in the second dual of $M(G)$, $M(G)^{**}$. By using [19, Corollary 1.2], we have $H^1(L^1(G),M(G))=0$. Then   by using proceeding corollary,  every $weak^*-to-weak$ continuous derivation from $L^1(G)^{**}$ into $M(G)$ is inner.\\
ii)  We know that  $c_0$ is  a $C^*-algebra$ and Since enery $C^*-algebra$ is weakly amenable, it follows that $c_0$ is weakly amenable. Then by using proceeding theorem, every $weak^*-to-weak$ continuous derivation from $\ell_\infty$ into $\ell^1$ is inner.\\\\

\noindent{\it{\bf Remark 4-5.}} Now let $B$ be a Banach $A-bimodule$ and for every derivation $D:A^{**}\rightarrow B$, we have $D(A^{**})\subseteq \overline{D\mid_A(A)}^w$. Assume that $D\in Z^1(A^{**},B)$. Take $\widetilde{D}=D\mid_A$. Then $\widetilde{D}\in Z^1(A,B)$. If $H^1(A,B)=0$, then there is a $b\in B$ such that $\widetilde{D}=\delta_b$.  Since $D(A^{**})\subseteq \overline{D\mid_A(A)}^w$, for every $a^{\prime\prime}\in A^{**}$, there is $(a_{\alpha})_{\alpha}\subseteq A$ such that
 $\widetilde{D} (a_{\alpha}) \stackrel{w} {\rightarrow}D(a^{\prime\prime})$. Thus we have the following equality.
 $$D(a^{\prime\prime})=w-\lim_\alpha\widetilde{D}(a_\alpha )=w-\lim_\alpha\delta_b(a_\alpha )=w-\lim_\alpha (a_\alpha b-  ba_\alpha)~in~B.$$   Thus, for every derivation $D:A^{**}\rightarrow B$, by assumption $D(A^{**})\subseteq \overline{D\mid_A(A)}^w$ and $H^1(A,B)=0$, we can write $D(a^{\prime\prime})=\lim_\alpha (a_\alpha b-  ba_\alpha)~in~B.$\\\\

\noindent{\it{\bf Theorem 4-6.}}  Let $B$ be a Banach $A-bimodule$ and $A$ has a $LBAI$. Suppose that $AB^{**},~B^{**}A\subseteq B$ and  every derivation $A^{**}$ into  $B^*$ is $weak^*-to-weak^*$ continuous. Then if $H^1(A,B^*)=0$, it follows that  $H^1(A^{**},B^{***})=0$.
\begin{proof} Take $B^{***}=B^*\oplus B^\perp$ where $B^\perp=\{b^{\prime\prime\prime}\in B^{***}:~b^{\prime\prime\prime}\mid_B=0\}$. Then we have
$$H^1(A^{**},B^{***})=H^1(A^{**},B^*)\oplus H^1(A^{**},B^\perp).$$
Since $H^1(A,B^*)=0$, it is similar to Theorem 2.17  that $H^1(A^{**},B^*)=0$. It suffices for the result to show that $H^1(A^{**},B^\perp)=0$.
 Without loss generality, let $(e_{\alpha})_{\alpha}\subseteq A$ be a $LBAI$ for $A$ such that  $e_{\alpha} \stackrel{w^*} {\rightarrow}e^{\prime\prime}$ in $A^{**}$ where $e^{\prime\prime}$ is a left unit for  $A^{**}$ with respect to the first Arens product. Let $a^{\prime\prime}\in A^{**}$ and suppose that $(a_{\beta})_{\beta}\subseteq A$ such that
 $a_{\beta} \stackrel{w^*} {\rightarrow}a^{\prime\prime}$ in $A^{**}$. Then, if we take $D\in Z^1(A^{**},B^\perp)$, for every $b^{\prime\prime}\in B^{**}$, by using the $weak^*-to-weak^*$ continuity  of $D$, we have
$$\langle D(a^{\prime\prime}), b^{\prime\prime}\rangle= \langle D(e^{\prime\prime}a^{\prime\prime}), b^{\prime\prime}\rangle=
\lim_\alpha\lim_\beta \langle (D(e_\alpha a_\beta),b^{\prime\prime}\rangle$$$$
=\lim_\alpha\lim_\beta \langle (D(e_\alpha) a_\beta +e_\alpha D(a_\beta)),b^{\prime\prime}\rangle
=\lim_\alpha\lim_\beta \langle D(e_\alpha) a_\beta ,b^{\prime\prime}\rangle$$$$+\lim_\alpha\lim_\beta \langle e_\alpha D(a_\beta),b^{\prime\prime}\rangle
=\lim_\alpha\lim_\beta \langle D(e_\alpha), a_\beta b^{\prime\prime}\rangle$$$$+\lim_\alpha\lim_\beta \langle  D(a_\beta),b^{\prime\prime}e_\alpha \rangle=0.$$
It follows that $D=0$, and so the result is hold.\\ \end{proof}

\noindent{\it{\bf Corollary 4-7.}} Assume that $A$ is a Banach algebra with $LBAI$. Suppose that $A$ is two-sided ideal in $A^{**}$  and every derivation  $D:A^{**}\rightarrow A^{***}$ is $weak^*-to- weak^*$ continuous. Then if $A$ is weakly amenable, it follows that $A^{**}$ is weakly amenable.\\\\

Let $n\geq 0$ and suppose that $A$ is $(2n+1)-weakly$ amenable.  By conditions of the proceeding corollary and by using [6, Corollary 1.14], we conclude that $A^{**}$ is weakly amenable.\\
Assume that $G$ is a locally compact group. We know that $L^1(G)$ is weakly amenable Banach algebra, see [16]. So by using proceeding corollary,  every $weak^*-to-weak^*$ continuous derivation from $L^1(G)^{**}$ into $L^1(G)^{***}$ is inner.\\

Let $A$ be a Banach Algebra. A dual Banach $A-bimodule$ $B$ is called normal if, for each $b\in B$ the map $a\rightarrow ax$ and $a\rightarrow xa$ from $A$ into $B$ is $weak^*-to-weak^*$ continuous.\\
A dual Banach algebra $A$ is Connes-amenable if, for every normal, dual Banach $A-bimodule$ $B$, every $weak^*-to-weak^*$ continuous derivation $D\in Z^1(A,B)$ is inner. Then we write $H^1_{w^*}(A,B)=0$.\\
A Banach algebra $A$ is called super-amenable if $H^1(A,B)=0$ for every Banach $A-bimodule$ $B$. It is clear that if $A$ is super-amenable, then $A$ is  amenable.\\
If $B$ is a Banach algebra and it is consisting a closed subalgebra $A$ such that $A$ is weakly dense in $B$. Then, it is easy to show that $A$ is a super-amenable if and only if $B$ is a super-amenable.\\
In the following, we will study some problems on the Connes-amenability of Banach algebras. In the  Theorem 4-3, with some conditions, we showed that $H^1(A^{**},B^{**})=0$. Now,  in the following theorem, by using some conditions and super-amenability of Banach algebra $A$, we conclude that $H^1_{w^*}(A^{**},B^{**})=0$.\\
By Theorem 4.4.8 from [25], we know that if $A$ is Arens regular Banach algebra and $A$ is an ideal in $A^{**}$, then $A$ is amenable if and only if $A^{**}$ is Connes-amenable. Now in the following by using of super-amenability of Banach algebra $A$, we show that $A^{**}$ is Connes-amenable.\\
In the Theorem 4-3, for a Banach $A-bimodule$ $B$, with some conditions, we showed that $H^1(A^{**},B^{**})=0$. In the following, by using some new conditions, we show that $H_{w^*}^1(A^{**},B^{**})=0$ and as results in group algebra, for a  amenable compact group $G$, we show that $H^1_{w^*}(L^1(G)^{**},M(G)^{***})=0.$\\\\

\noindent{\it{\bf Theorem 4-8.}} Suppose that $A$ is supper-amenable  Banach algebra. Then we have the following statements.
\begin{enumerate}
\item $A^{**}$ is Connes-amenable Banach algebra with respect to the first Arens product.
\item  If $A$ has a $LBAI$ and for every Banach $A-bimodule$ $B$ with predual $X$, we have $AB^*,~B^*A\subseteq X$, then $H^1_{w^*}(A^{**},B^{**})=0$.
\end{enumerate}
\begin{proof}
\begin{enumerate}
\item Let  $B$ be a normal dual Banach $A-bimodule$ $B$ and let  be $D\in Z^1(A^{**},B^{**})$ be $weak^*-to-weak^*$ continuous derivation. Suppose that $\widetilde{D}=D\mid_A$. Then we have  $\widetilde{D}\in Z^1(A,B)$. Since $A$ is super-amenable, there is element  $b\in B$ such that $\widetilde{D}=\delta_b$. Suppose that $a^{\prime\prime}\in A^{**}$ and $(a_{\alpha})_{\alpha}\subseteq A$  such that  $a_{\alpha} \stackrel{w^*} {\rightarrow}a^{\prime\prime}$ in $A^{**}$. Since $\widetilde{D}$ is $weak^*-to-weak^*$ continuous derivation, we have
$$D(a^{\prime\prime})=D(w^*-\lim_\alpha a_\alpha)=w^*-\lim_\alpha D(a_\alpha)=w^*-\lim_\alpha \widetilde{D}(a_\alpha)$$$$=w^*-\lim_\alpha \delta_b(a_\alpha)=w^*-\lim_\alpha (a_\alpha b-ba_\alpha )=a^{\prime\prime}b-ba^{\prime\prime}.$$
It follows that $H^1_{w^*}(A^{**},B)=0$.

\item By using part (1), proof is similar to Theorem 4-3 (2).\\
\end{enumerate}
\end{proof}

\noindent{\it{\bf Theorem 4-9.}} Suppose that $A$ is  an amenable  Banach algebra with $LBAI$. If for every Banach $A-bimodule$ $B$, we have $AB^{**},~B^{**}A\subseteq B$, then $H^1_{w^*}(A^{**},B^{***})=0$.
\begin{proof} Proof is similar to the proceeding theorem.\\\end{proof}

\noindent{\it{\bf Example 4-10.}} Suppose that $G$ is a compact group. Then $L^1(G)$ is an ideal in $M(G)^{**}$. If $G$ is amenable, then by using proceeding theorem, we conclude the following equality.
$$H^1_{w^*}(L^1(G)^{**},M(G)^{***})=0.$$\\

\noindent{\it{\bf Corollary 4-11.}} Assume that $A$ is a weakly amenable Banach algebra with $LBAI$. Then if $A$ is an ideal in $A^{**}$, it follows that
$$H^1_{w^*}(A^{**},A^{***})=0.$$\\

\noindent{\it{\bf Example 4-12.}} Assume that $G$ is a compact group. Then we know that $L^1(G)$ has a bounded approximate identity [=$BAI$] and $L^1(G)$ is two-sided ideal in $L^1(G)^{**}$.  We know that $L^1(G)$ is weakly amenable. Hence it is clear $$H^1_{w^*}(L^\infty(G)^*,L^\infty(G)^{**})=0.$$\\

\begin{center}
\section{ \bf A representations of derivations and Arens regularity od Banach algebras  }
\end{center}

Let $B$ be a Banach $A-bimodule$. Then for every $b\in B$, we define
$$L_{b}(a)=ba~~~and~~~~
R_{b}(a)=ab,$$
for every $a\in A$. These are the operation of left and right multiplication by $b$ on $A$.\\
In the following, for super-amenable Banach algebra $A$ and a Banach $A-bimodule$ $C$, we have represented for every derivation $D$ from $A$ into $C$, and as a result in group algebras, for amenable locally compact group $G$, we show that every derivation $D$ from $L^1(G)$ into $L^\infty (G)$ is in the form $D=L_f$ where $f\in L^\infty (G)$.\\
In the following by using the super-amenability of Banach algebra $A$, we give a representation  for $Z^1(A,C)$ where $C$ is a Banach $A-bimodule$. As corollary, For an amenable locally compact group $G$. We show that  there is a Banach $L^1(G)-bimodule$ such as $(L^\infty(G),.)$ such that $Z^1(L^1(G),L^\infty(G))=\{L_{f}:~f\in L^\infty(G)\}.$\\\\

\noindent{\it{\bf Theorem 5-1.}} Assume that $A$ is a super-amenable Banach algebra and $B$ is a Banach $A-bimodule$ such that $B$ factors on the left [resp. right]. Then there is a Banach $A-bimodule$ $(C,.)$ such that $C=B$ and
$Z^1(A,C)=\{L_{D^{\prime\prime}(e^{\prime\prime})}:~D\in Z^1(A,C)\}~~[resp.~~Z^1(A,C)=\{R_{D^{\prime\prime}(e^{\prime\prime})}:~D\in Z^1(A,C)\}]$ where $e^{\prime\prime}$ is a right [resp. left] identity for $A^{**}$.\\
\begin{proof} Certainly, every super-amenable Banach algebra is amenable. So, by using [25, Proposition 2.2.1], $A$ has a $BAI$ such as $(e_{\alpha})_{\alpha}$. Since $BA=B$, for every $b\in B$, there are $y\in B$ and $a\in A$ such that $b=ya$. Then we have
$$\lim_\alpha b e_{\alpha}=\lim_\alpha (ya) e_{\alpha}=\lim_\alpha y(a e_{\alpha})=ya=b.$$
It follows that $B$ has $RBAI$ as $(e_{\alpha})_{\alpha}\subseteq A$. Without loss generality, let $e^{\prime\prime}$ be a right unit for $A^{**}$ such that $e_{\alpha} \stackrel{w^*} {\rightarrow}e^{\prime\prime}$ in $A^{**}$.\\
Take $C=B$ and for every $a\in A$ and $x\in C$, we define $a.x=0$ and $x.a=xa$. It is clear that $(C,.)$ is a Banach $A-bimodule$. Suppose that $D\in Z^1(A,C)$. Then there is a $c\in C$ such that $D=\delta_c$. Then for every $a\in A$, we have
$$D(a)=\delta_c (a)=a.c-c.a=-ca.$$
Suppose that $x\in C$ and $x^\prime\in C^*$. Since $C.A=C$, there is $t\in C$ and $s\in A$ such that $x=t.s=ts$. Then we have
$$\langle x^\prime, xe_\alpha\rangle=\langle x^\prime, tse_\alpha\rangle=\langle x^\prime t,se_\alpha\rangle\rightarrow
\langle x^\prime t,s\rangle=\langle x^\prime, x\rangle.$$
It follows that $xe_{\alpha} \stackrel{w} {\rightarrow}x$ in $C$. Since $D^{\prime\prime}$ is a $weak^*-to-weak^*$ continuous derivation, we have
$$D^{\prime\prime}(e^{\prime\prime})=D^{\prime\prime}(w^*-\lim_\alpha e_{\alpha})=w^*-\lim_\alpha D^{\prime\prime}(e_{\alpha})=w-\lim_\alpha D(e_{\alpha})$$$$=w-\lim_\alpha (-ce_{\alpha})=-c.$$
Thus we conclude that $D(a)=D^{\prime\prime}(e^{\prime\prime})a$ for all $a\in A$. It follows that $$D=L_{D^{\prime\prime}(e^{\prime\prime})}.$$
On the other hand, since for every derivation $D\in Z^1(A,C)$,  $L_{D^{\prime\prime}(e^{\prime\prime})}\in
Z^1(A,C)$, the result  hold.\\\end{proof}

\noindent{\it{\bf Corollary 5-2.}} Suppose that $A$ is  an amenable  Banach algebra and $B$ is a Banach $A-bimodule$ such that $AB^*=B^*$ [resp. $B^*A=B^*$]. Then there is a Banach $A-bimodule$ $(C,.)$ such that $C=B$ and
$Z^1(A,C^*)=\{L_{D^{\prime\prime}(e^{\prime\prime})}:~D\in Z^1(A,C^*)\}~~[resp.~~Z^1(A,C^*)=\{R_{D^{\prime\prime}(e^{\prime\prime})}:~D\in Z^1(A,C^*)\}]$ where $e^{\prime\prime}$ is a right [resp. left] identity for $A^{**}$.\\\\

\noindent{\it{\bf Corollary 5-3.}} Suppose that $A$ is  a super-amenable  Banach algebra and $B$ is a Banach $A-bimodule$. Then there is a Banach $A-bimodule$ $(C,.)$ such that $C=B$ and
$Z^1(A,C)=\{L_{c}:~c\in C\}~~[resp.~~Z^1(A,C)=\{R_{c}:~c\in C\}]$.\\\\

\noindent{\it{\bf Example 5-4.}}
\begin{enumerate}
\item  Let $G$ be an amenable locally compact group. Then there is a Banach $L^1(G)-bimodule$ such as $(L^\infty(G),.)$ such that $Z^1(L^1(G),L^\infty(G))=\{L_{f}:~f\in L^\infty(G)\}.$
\item  Let $G$ be locally compact group and $H^1(M(G),L^1(G))=0$. Then there is a Banach $M(G)-bimodule$ $(L^1(G),.)$ such that $Z^1(M(G),L^1(G))\\=\{L_{D^{\prime\prime}(e^{\prime\prime})}:~D\in Z^1(M(G),L^1(G))\}~~[or~~Z^1(M(G),L^1(G))=\{R_{D^{\prime\prime}(e^{\prime\prime})}:~\\D\in Z^1(M(G),L^1(G))\}]$ where $e^{\prime\prime}$ is a left [or right] identity for $L^1(G)^{**}$.\\\\

\end{enumerate}
\noindent{\it{\bf Theorem 5-5.}}  Let $B$ be a Banach $A-bimodule$ and suppose that $D:A\rightarrow B^*$ is a derivation. If $D^{\prime\prime}:A^{**}\rightarrow B^{***}$ is a derivation and $B^*\subseteq D^{\prime\prime}(A^{**})$, then $Z^\ell_{A^{**}}(B^{**})=B^{**}$.
\begin{proof} Since $D^{\prime\prime}:A^{**}\rightarrow B^{***}$ is a derivation, by [20, Theorem 4.2],
$D^{\prime\prime}(A^{**})B^{**}\subseteq B^*$. Due to $B^*\subseteq D^{\prime\prime}(A^{**})$, we have $B^*B^{**}\subseteq B^*$. Let   $(a_{\alpha}^{\prime\prime})_{\alpha}\subseteq A^{**}$ such that  $a^{\prime\prime}_{\alpha} \stackrel{w^*} {\rightarrow}a^{\prime\prime}$ in $A^{**}$. Assume that $b^{\prime\prime}\in B^{**}$. Then for every $b^\prime\in B^*$, since $b^\prime b^{\prime\prime}\in B^*$, we have
$$\langle b^{\prime\prime}a_{\alpha}^{\prime\prime},b^\prime\rangle=\langle a_{\alpha}^{\prime\prime}b^\prime, b^{\prime\prime}\rangle\rightarrow
\langle a^{\prime\prime},b^\prime b^{\prime\prime}\rangle=\langle b^{\prime\prime}a^{\prime\prime},b^\prime \rangle.$$
Thus $b^{\prime\prime}a^{\prime\prime}_{\alpha} \stackrel{w^*} {\rightarrow}b^{\prime\prime}a^{\prime\prime}$ in $B^{**}$, and so $b^{\prime\prime}\in Z^\ell_{A^{**}}(B^{**})$.\\ \end{proof}

\noindent{\it{\bf Corollary 5-6.}} Assume that $A$ is a Banach algebra and $D:A\rightarrow A^*$ is a derivation such that $A^*\subseteq D^{\prime\prime}(A^{**})$. Then if the linear mapping $D^{\prime\prime}:A^{**}\rightarrow A^{***}$
is a derivation, it follows that  $A$ is Arens regular.\\\\

\noindent{\it{\bf Lemma 5-7.}} Let $B$ be a Banach left $A-module$ and $B^{**}$ has a $LBAI$ with respect to $A^{**}$. Then $B^{**}$ has a left unit with respect to $A^{**}$.
\begin{proof}
Assume that $(e^{\prime\prime}_{\alpha})_{\alpha}\subseteq A^{**}$ is a $LBAI$ for $B^{**}$. By passing to a subnet, we may suppose that there is  $e^{\prime\prime}\in A^{**}$ such that $e^{\prime\prime}_{\alpha} \stackrel{w^*} {\rightarrow}e^{\prime\prime}$ in $A^{**}$. Then for every $b^{\prime\prime}\in B^{**}$ and $b^\prime\in B^*$, we have

$$\langle \pi_\ell^{***}(e^{\prime\prime},b^{\prime\prime}),b^\prime\rangle=
\langle e^{\prime\prime},\pi_\ell^{**}(b^{\prime\prime},b^\prime)\rangle=
\lim_\alpha \langle e_\alpha^{\prime\prime},\pi_\ell^{**}(b^{\prime\prime},b^\prime)\rangle$$$$=
\lim_\alpha \langle \pi_\ell^{***}(e_\alpha^{\prime\prime},b^{\prime\prime}),b^\prime\rangle=
\langle b^{\prime\prime},b^\prime\rangle.$$
It follows that $\pi_\ell^{***}(e^{\prime\prime},b^{\prime\prime})=b^{\prime\prime}$.\\\end{proof}

In the proceeding theorem, if we take $B=A$, then we obtain Theorem 1.1, from [12].\\\\

\noindent{\it{\bf Theorem 5-8.}} Let $A$ be a left strongly Arens irregular and suppose that $A^{**}$ is an amenable Banach algebra. Then we have the following assertions.
\begin{enumerate}
\item $A$ has an identity.

\item If $A$ is a dual Banach algebra, then $A$ is reflexive.
\end{enumerate}
\begin{proof}\begin{enumerate}
\item Since $A^{**}$ is amenable, it has a $BAI$.  By using  proceeding Lemma, $A^{**}$ has an identity $e^{\prime\prime}$. So, the mapping $x^{\prime\prime}\rightarrow  e^{\prime\prime} x^{\prime\prime}=x^{\prime\prime}$ is $weak^*-to-weak^*$ continuous from $A^{**}$ into $A^{**}$. It follows that $e^{\prime\prime}\in Z_1(A^{**})=A$. Consequently, $A$ has an identity.

\item Assume that $E$ is predual Banach algebra for $A$. Then we have $A^{**}=A\oplus E^\bot$. Since $A^{**}$ is amenable, by Theorem 2.3 from [11],  $A$ is amenable, and so $E^\bot$ is amenable. Thus $E^\bot$ has a $BAI$ such as $(e^{\prime\prime}_\alpha)_\alpha\subseteq E^\bot$. Since $E^\bot$ is a closed and $weak^*-closed$ subspace of $A^{**}$, without loss generality,  there is $e^{\prime\prime}\in E^\bot$ such that
   $$e^{\prime\prime}_\alpha\stackrel{w^*} {\rightarrow}e^{\prime\prime}~~and~~e^{\prime\prime}_\alpha\stackrel{\parallel~~\parallel} {\rightarrow}e^{\prime\prime}$$
 Then   $e^{\prime\prime}$ is a left identity for  $E^\bot$. On the other hand, for every $x^{\prime\prime}\in E^\bot$, since $E^\bot$ is an ideal in $A^{**}$, we have $x^{\prime\prime}e^{\prime\prime}\in E^\bot$. Thus for every $a^\prime\in A^*$, we have
 $$\langle x^{\prime\prime}e^{\prime\prime},a^\prime\rangle=\lim_\alpha
 \langle (x^{\prime\prime}e^{\prime\prime})e^{\prime\prime}_\alpha,a^\prime\rangle=\lim_\alpha
 \langle x^{\prime\prime}(e^{\prime\prime}e^{\prime\prime}_\alpha),a^\prime\rangle
 $$$$=\lim_\alpha
 \langle x^{\prime\prime}e^{\prime\prime}_\alpha,a^\prime\rangle= \langle x^{\prime\prime},a^\prime\rangle.$$
 It follows that $x^{\prime\prime}e^{\prime\prime}=x^{\prime\prime}$, and so  $e^{\prime\prime}$ is a right identity for $E^\bot$. Consequently, $e^{\prime\prime}$ is a two-sided identity for $E^\bot$.
Now, let $a^{\prime\prime}\in A^{**}$. Then we have the following equalities.
$$e^{\prime\prime}a^{\prime\prime}= (e^{\prime\prime}a^{\prime\prime})e^{\prime\prime}=
 e^{\prime\prime}(a^{\prime\prime}e^{\prime\prime})=a^{\prime\prime}e^{\prime\prime}.$$
Thus we have $e^{\prime\prime}\in Z_1(A^{**})=A$. It follows that  $e^{\prime\prime}=0$, and so $E^\bot=0$. Consequently, we have $A^{**}=A$.\\\end{enumerate}\end{proof}

\noindent{\it{\bf Example 5-9.}}  Let $G$ be a locally compact group. Then if $L^1(G)^{**}$ or $M(G)^{**}$ are amenable, by using proceeding theorem part (1) and (2), respectively,  we conclude that $G$ is finite group, see [12].\\\\

\noindent{\it{\bf Problems.}} \\
i) By notice to Example 2-13, for locally compact group $G$, if $e^{\prime\prime}\in L^1(G)^{**}$ is a mixed unit for $L^1(G)^{**}$, then find the following statements\\
$Z^\ell_{e^{\prime\prime}}(L^1(G)^{**})=?$ and $Z^\ell_{e^{\prime\prime}}(M(G)^{**})=?$\\
ii) Dose the following assertions hold?
\begin{enumerate}
\item For compact group $G$, we have $H^1(L^1(G)^{**},M(G))=0.$
\item  $H^1(\ell^\infty,\ell^1)=0.$\\\\
\end{enumerate}

\bibliographystyle{amsplain}

\noindent Department of Mathematics, University of Mohghegh Ardabili, Ardabil, Iran\\
{\it Email address:} haghnejad@aut.ac.ir

\end{document}